\documentclass[12pt,a4paper]{amsart}

\usepackage{bbm}
\usepackage{lscape}
\usepackage{tabularx,mdwtab}

\usepackage{hyperref}

 \def\mmat #1,#2,#3,#4,{\text{\small\arraycolsep=3pt $
\begin{pmatrix}#1&#2\\#3&#4\end{pmatrix}$}}

\newcommand{\un}{\underline{N}}

\newcommand {\fbj}{{\mathfrak{bj}}}

\newcommand {\fbr}{{\mathfrak{br}}}
\newcommand {\fbrj}{{\mathfrak{brj}}}

\newcommand {\fout}{{\mathfrak{out}}}
\newcommand {\fwk}{{\mathfrak{wk}}}

\newcommand {\fer}  {{\mathfrak{er}}}

\newcommand {\ffr}   {{\mathfrak{fr}}}

\newcommand {\fdy}   {{\mathfrak{dy}}}
\newcommand {\fby}   {{\mathfrak{by}}}
\newcommand {\fsby}   {{\mathfrak{sby}}}
\newcommand {\fmy}   {{\mathfrak{my}}}
\newcommand {\fsmy}   {{\mathfrak{smy}}}
\newcommand {\Lra}  {{\Longrightarrow}}

\newcommand{\del}{\partial}

\usepackage{DLdef1}

\begin{document}

\title{Non-degenerate invariant (super)symmetric bilinear forms on simple Lie (super)algebras}

\author{Sofiane Bouarroudj${}^a$, Andrey Krutov${}^{b,d}$,  \\
Dimitry Leites${}^{a,c}$, Irina
Shchepochkina${}^d$}

\address{${}^a$New York University Abu Dhabi,
Division of Science and Mathematics, P.O. Box 129188, United Arab
Emirates; sofiane.bouarroudj@nyu.edu\\
$^b$Institute of Mathematics, Polish Academy of Sciences, ul. \'{S}niadeckich 8, 00-656 Warszawa, Poland;
a.krutov@impan.pl\\
${}^{c}$Department of Mathematics, Stockholm University, SE-106 91
Stockholm, Sweden;
mleites@math.su.se; dl146@nyu.edu\\
${}^d$Inde\-pendent University of Moscow, Bolshoj Vlasievsky per, dom
11, RU-119 002 Moscow, Russia; irina@mccme.ru}


\begin{abstract}
We review the list of non-degenerate invariant (super)symmetric bilinear forms (briefly: NIS) on the following simple (relatives of) Lie (super)algebras: (a) with symmetrizable Cartan matrix of any growth, (b) with non-symmetrizable Cartan matrix of polynomial growth, (c) Lie (super)algebras of vector fields with polynomial coefficients, (d) stringy a.k.a. superconformal superalgebras, (e)  queerifications of simple restricted Lie algebras.

Over algebraically closed fields of positive characteristic, we establish when the deform (i.e., the result of deformation) of the known finite-dimensional simple Lie (super)algebra has a NIS. Amazingly, in most of the cases considered, if the Lie (super)algebra has a NIS, its deform has a NIS with the same Gram matrix after an identification of bases of the initial and deformed algebras. We do not consider odd parameters of deformations.

Closely related with simple Lie (super)algebras with NIS is the notion of doubly extended Lie (super)algebras of which affine Kac--Moody (super)algebras are the most known examples.
\end{abstract}

\thanks{S.B. and A.K. were partly supported by the grant AD 065 NYUAD.  A.K. was partly supported by WCMCS post-doctoral fellowship. A part of this research was done while A.K. was visiting NYUAD; the
financial support and warm atmosphere of this institute are gratefully
acknowledged. We are thankful to J.~Bernstein, P.~Grozman, S.~Skryabin,
 P.~Zusmanovich, and especially A.~Lebedev,  for help. For the possibility to conduct difficult computations of this research we are grateful to  M.~Al Barwani, Director of the High Performance Computing resources at New York University Abu Dhabi.}


\maketitle

\markboth{\itshape Sofiane Bouarroudj\textup{,} Andrey Krutov\textup{,} Dimitry Leites\textup{,} Irina Shchepochkina}{{\itshape Nondegenerate invariant (super)symmetric bilinear forms on simple Lie (super)algebras}}

\begin{quote}\rightline{To Alexandre Kirillov-p\`{e}re, our teacher}\end{quote}

\thispagestyle{empty}

\setcounter{tocdepth}{3}
\tableofcontents

\section{Introduction}

This is a sequel to our talk at the conference \lq\lq Representation Theory at the Crossroads of Modern Mathematics\rq\rq\
in honor of Alexandre Kirillov-p\`{e}re,
Reims, May 29 -- June 2, 2017. At the talk, we reported two results: that of \cite{BLLSq} (classification of the simple Lie superalgebras in characteristic 2 --- modulo classification of simple Lie algebras over the same field --- by means of the two methods: queerification and \lq\lq method 2\rq\rq) and  that of  \cite{KrLe} vital for \lq\lq method 2\rq\rq.

We needed
derivations and central extensions of simple Lie (super)algebras in the approach to the classification of simple $\Zee$-graded Lie (super)algebras of depth 1 we used over $\Cee$, see \cite{LSh}; we hope to adjust the same approach to classification of simple vectorial Lie (super)algebras for $p=3$ and 2 (for $p>3$, the classification is obtained; for a summary, see \cite{S}). For certain types of simple Lie (super)algebras their derivations and central extensions  are described in \cite{BGLL},
where several results are interpreted using a recently distinguished notion: \textbf{double extension} of Lie (super)algebras, see \S~\ref{DE}. One of the three ingredients that define a \textit{double extension} of a given Lie (super)algebra $\fh$ is a non-degenerate invariant (super)symmetric bilinear form (briefly: NIS) on $\fh$.

Here, we consider two types of $\Zee$-graded simple (or close to simple) Lie (super)algebras for which classification is obtained or at least conjectured: finite-dimensional  and  of polynomial growth; we mention several other types as well. We recall the known results (which of the algebras of these two types have a NIS) and, for the cases where the \textbf{true deforms} (the results of deformations, neither trivial, nor semitrivial;  for a characterization of semitrivial deforms in some cases, see \cite{BLLS} and especially \cite{BGLd}) of these algebras are classified, we investigate if the deform has a NIS. (We also have to mention several non-simple but natural examples of Lie (super)algebras, e.g., the ones with non-invertible Cartan matrix and results of (generalized) Cartan prolongations, see \cite{Shch}.)

\ssec{General comments}
Over fields $\Kee$ of characteristic $p>0$, simple Lie (super)algebras have overwhelmingly many deformations; for $p$ small (3 or, even worse, 2), the number of deformations becomes appalling.  We are interested, however, in \textbf{deforms}, i.e., the \textit{results} of deformations, rather than in deformations themselves, besides, some deforms, although correspond to integrable, and not cohomologous to zero, cocycles are \textbf{semi-trivial}, i.e., the corresponding deforms are isomorphic to the initial (super)algebra; only few deforms are \textbf{true}, see \cite{BLW, BGLd}.

On the last page of \cite{Zus}, there are given conditions  for the deform of a Lie algebra $\fg$ with the new bracket $[x,y]_t=\sum \varphi_i(x,y)t^i$, where $\varphi_0(x,y)=[x, y]$, to have a NIS $(x,y)_t=\sum (x,y)_jt^j$, provided $(\cdot ,\cdot)_0$ is a NIS on $\fg$. Namely, for any $x,y,z\in \fg$ we should have
\[
\mathop{\sum}\limits_{i+j=n, \ 0\leq i, j}(\varphi_i(x,y), z)_j+(x, \varphi_i(y, z))_j=0\text{~~ for $n=0, 1, \dots$}.
\]
In our computations (proofs of Claims with the aid of the \textit{SuperLie} package \cite{Gr}) we checked when these conditions, and their (not obvious) superizations are fulfilled. We observed the following very interesting fact we can not explain.

\sssec{Fact}\label{Fact} \textit{In all cases for $p\neq2$ we know, except the serial deforms $\fsvect_{\exp_i}$ and $\fsvect_{1+\bar u}^{(1)}$, see table~\eqref{svectH}, if a given simple Lie (super)algebra has a NIS, then its deform also has a NIS.}

The results of computer experiments (the two series of exceptions from Fact, see Conjecture ~\ref{ClSvect}) trouble us, but these exceptions are  also a fact (and \lq\lq fact are stubborn things").

In  \cite[Proposition~5.10 (Proposition~4.4.1 in the arXiv version)]{PU},  it is proved, using simple arguments, that \lq\lq a deformation of a finite-dimensional complex Lie algebra with a NIS, is equivalent to a deformation with the form unchanged..." The statement in \cite{PU} is not astonishing: over a quadratically closed fields of characteristic $\neq 2$ any NIS has the unit Gram matrix in some basis. If we deform this form, then at first it remains non-degenerate, and it is conceivable to change a basis to keep the same Gram matrix. But why is it possible to deform the form preserving its invariance? So Fact~\ref{Fact} is remarkable, whereas exceptions from it should have been expected.

\ssec{The types of Lie (super)algebras we consider}  (A) \underline{Over $\Cee$}.

If $A$  is symmetrizable and invertible, then the Lie (super)algebra $\fg(A)$ has a NIS, see~\eqref{recipe}, provided either $p\neq 2$, or $p=2$ and $\fg_{2\alpha}=0$ for any root $\alpha$.

If $A$  is non-symmetrizable corresponding to the Lie superalgebra $\fg(A)$ of polynomial growth, there are two types of such Lie superalgebras; on $\fpsq^{\ell(2)}(n)$, there is an odd NIS for $\fg(A)$ of one type, and no NIS for $\fg(A)$ of the other type, see \S~7. 

\underline{Over algebraically closed fields $\Kee$ of characteristic $p>0$},  finite-dimensional Lie algebras and Lie (super)algebras $\fg(A)$ with indecomposable Cartan matrix $A$ have been classified, see \cite{BGL}.

We investigate the existence of NIS on \textit{all} true deforms of Lie algebras of the form $\fg(A)$, and on their simple relatives. We do not know how many non-isomorphic Lie algebras are there among the deforms of a given $\fg$. To answer this question, one should apply (for example) the technique of Kuznetsov and Chebochko, cf. \cite{KuCh, ChKu}. In the absence of the answer, we chose a basis in the cohomology space, and for every element of the basis, we tested a cocycle representing it.

(B) To Dzhumadildaev's announcement \cite{Dz1}  on classification of NISes on simple $\Zee$-graded vectorial Lie algebras in characteristic $p>0$, reproduced (with explicit proof  added) in the book by  Strade and Farnsteiner \cite{SF}\footnote{In \cite{Fa,SF}, the \textbf{invariant bilinear form on} the Lie algebra $\fg$ is referred to as \textit{associative form of} $\fg$.}   we add

(Ba) the examples of simple vectorial Lie (super)algebras for $p=5$ and $3$ found after \cite{Dz1, SF} were published, and

(Bb) deforms of simple finite-dimensional Lie (super)algebras $\fg(A)$ with indecomposable Cartan matrix $A$ or simple subquotients (\lq\lq relatives") thereof; (for the cases these deforms  have no Cartan matrix) classified
for $p>0$, see \cite{BGL, BGLd}. 

(C) We extend these investigations of cases (A) and (B) to simple vectorial Lie \textbf{super}algebras, first of all over $\Cee$, where the classification is obtained (see \cite{Sh5, Sh14, LSh} and \cite{K,K10, CaKa} and references therein) and also over $\Kee$ for $p>2$ as far as a conjecture and partial classification results (see \cite{BGLLS, BGLLS1, BLLS}) enable us to embrace the known examples and correct some claims in \cite{NS}.

(D) We also consider the simple Lie superalgebras with an odd supersymmetric invariant bilinear form, and their non-simple relatives, both finite-dimensional and of polynomial growth.

\textbf{We do not consider odd parameters of deformations}.

\ssec{NIS: selected applications}\label{ssPub} In (quantum)
inverse scattering method for solving partial differential equations: for a general review,
see~\cite[\S 1.1, \S 4.1]{FT}, \cite{BS}; for a
discussion in the context of BV-quantization, see~\cite{KK}.

It is known (due to P. Etingof, I. Losev) that the (super)algebra of observables (deformed Fock space) of the Calogero models based on the root system is simple for almost all values of the coupling constants. Konstein et al. found exact values of coupling constants for which the bilinear forms induced by (super)traces become degenerate, and the (super)algebra of observables acquires ideals, see \cite{KS,KT}.


\section{Double extension: an interesting notion}\label{DE}

Well-known examples of \textit{indecomposable  double extensions} to be defined shortly are affine Kac--Moody algebras $\fg(A)$ with indecomposable Cartan matrix $A$ over $\Cee$, see \cite{Kb}; $\fgl(pn)$ 
and  $\fgl(n|n+kp)$ for any $k\in\Nee$ in characteristic $p>0$; Lie superalgebras $\fgl(n|n)$, $\fq(n)$, and $\fpo(0|m)$ for any $p$. \lq\lq Pendant plus de quarante ans, nous avons parl\'e de prose --- et ce n'est pas connu!"

 \ssec{Definition of double extension}\label{DD} In this Section, we explain how the construction of \textit{double extension} works for \textbf{Lie algebras over fields of any characteristic, and Lie superal\-geb\-ras over fields of characteristic $p\neq 2$}. For the (rather non-trivial) case $p= 2$, see \cite{BeBou}.

Observe that for Lie algebras over fields $\Kee$ of characteristic $p>0$ and Lie superalgebras for any $p$, an analog of well-known NIS, called the \textit{Killing form}, should be sought in projective representa\-tions.  In other words, the said analog should be related to a non-trivial central extension, see \cite{Kapp, GP}. It only remained to add a(n outer) derivation to get even nicer object. The double extension requires one more ingredient: a NIS.

Given a Lie (super)algebra $\fh$, its \textit{double extension} $\fg$ simultaneously involves three ingredients:

1) a central extension $\fh_c$ of $\fh$ with the center spanned by $c$, so $\fh\simeq\fh_c/\Kee c$,

2) a derivation $D$ of $\fh_c$ such that
$\fg\simeq\fh_c\ltimes \Kee D$, a semidirect sum,

3)  an $\fh$-invariant NIS $B_\fh$ on $\fh$. Observe that, in super setting, $B_\fh$ can be odd.

Then, under certain conditions, there is a NIS $B$ on $\fg$, which extends $B_\fh$, and the interesting cases are the ones where the Lie algebra $\fg$ is not a direct sum of ideals $\fh$ and $\Kee c\oplus \Kee D$;  we call such direct sums \textit{decomposable} double extensions; it is called \textit{reducible} in \cite{BeBou}.


\sssec{Lemma (On a central extension)}\textup{(Lemma 3.6, page 73 in \cite{BB})} \label{L1}
\textit{Let $\fh$ be a  Lie (super)algebra over a field $\Kee$, let $B_\fh$ be   an $\fh$-invariant NIS  on $\fh$, and $D\in\fder\ \fh$ a derivation  such that $B_\fh$ is $D$-invariant, i.e.,
\be\label{(1)}
B_\fh(Da,b)+(-1)^{p(a)p(D)}B_\fh(a,Db)=0\text{~~ for any $a,b\in \fh$. } \ee
Then the bilinear form $\omega(a,b):=B_\fh(Da,b)$ is a $2$-cocycle of the Lie (super)algebra $\fh$.}

Thus, in the assumptions of Lemma~\ref{L1}, we can construct a central extension $\fh_\omega$ of $\fh$ given by cocycle $\omega$  so that $\fh_\omega/\Kee c\simeq \fh$.

Let us find out what the conditions for $D$-invariance of  the cocycle $\omega$ are, i.e., when the operator $d\in\fgl(\fh_\omega)$, such that $d(c)=0$ and $d|_\fh=D$, is  a~derivation of the Lie (super)algebra $\fh_\omega$. For this, we have to verify that
\be\label{(5)}
\omega(Dx,y)+(-1)^{p(x)p(D)}\omega(x,Dy)=0 \text{ for all } x,y\in \fh.
\ee
We have
$$
\begin{array}{l}
\omega(Dx,y)+(-1)^{p(x)p(D)}\omega(x,Dy)=B_\fh(D^2x,y)+(-1)^{p(x)p(D)}B_\fh(Dx,Dy)=\\
B_\fh(D^2x,y)-(-1)^{p(x)p(D)+p(D)(p(D)+p(x))}
B_\fh(D^2x,y)=
(1-(-1)^{p(D)^2})B_\fh(D^2x,y).
\end{array}
$$
Since the form $B_\fh$ is nondegenerate on $\fh$, it follows that equality \eqref{(5)} holds for any $x,y\in\fh$ if and only if the operator $D$  is either even, or odd and such that
$D^2=0$. (Recall again that we do not consider Lie superalgebras for $p=2$.)

Thus, we have a Lie (super)algebra $\fh_\omega$ and its derivation $d$. Generally, given an arbitrary Lie superalgebra $\fh$ and its even derivation $D$, we can always construct a semidirect sum $\fh\ltimes \Kee D$. If, however, $D$ is odd, to construct a semidirect sum $\fh\ltimes \Kee D$, we must require that $D^2$ is an inner derivation. Luckily, even a stronger condition is satisfied: $d^2=0$, and hence we can construct a semidirect sum $\fg=\fh_\omega\ltimes \Kee d$.

On the Lie (super)algebra $\fg$, define a (super)symmetric form $B$ by setting
\be\label{nis}
B|_{\fh}=B_\fh, \quad B(d,c)=1, \quad B(c,x)=B(d,x)=B(c,c)=B(d,d)=0 \text{ for all } x\in\fh.
\ee

\sssec{Lemma (On non-degenerate invariant symmetric forms)} \textup{(Theorem 1, page 68 in \cite{BB}; for Lie algebras: Exercise 2.10 in \cite{Kb})} \label{L2}
\textit{The form $B$ defined by \eqref{nis} is a NIS on $\fg$.}

Thus constructed Lie (super)algebra $\fg$ with NIS $B$ on it is called the \textit{double extension} or, for emphasis, \textit{D-extension}, of $\fh$.

\parbegin{Remark} If the derivation $D$ is inner, i.e., there exists an $x\in \fh$ such that $D=\ad_x$, then the operator $D-\ad_x$ vanishes identically. Replacing  $D$ by $D-\ad_x$, we see that the cocycle $\omega$ is also the zero one, and hence the Lie (super)algebra $\fg$ is a direct sum of its ideal $\fh$ and a~2-dimensional commutative ideal, on which we defined a nondegenerate symmetric bilinear form.
Evidently, this ---  \textit{decomposable} --- case is not interesting.
\end{Remark}

\ssec{Lie (super)algebra $\fg$ that can be a double extension of a~ Lie (super)algebra $\fh$}
Let $\fg$ be a Lie algebra over any field  $\Kee$ or a Lie superalgebra over a field $\Kee$ of characteristic $p\ne 2$; let $B$ be a non-degenerate invariant (super)symmetric bilinear form on  $\fg$, and $c\ne 0$ a central element of $\fg$.

The invariance
of the form $B$  implies that $B(c,[x,z])=0$ for any $x,z\in\fg$, i.e., the space $c^\perp$ contains the commutant $\fg^{(1)}=[\fg,\fg]$ of $\fg$,  and hence is an ideal. Since the form $B$ is nondegenerate, the codimension of this ideal is equal to 1.

If $B(c,c)\ne 0$, then the Lie (super)algebra $\fg$ is just a direct sum $\fg=\Kee c\oplus c^\perp$. This case is not interesting, i.e., \textit{decomposable}.

Moreover, even if $B(c,c)= 0$, but $c\notin \fg^{(1)}$, then any subspace $V\subset\fg$, complementing $\Kee c$ and containing $\fg^{(1)}$, is an ideal and the Lie (super)algebra $\fg$ is a direct sum of ideals: $\fg=V\oplus \Kee c$. On the ideal $V$, the form $B$ must be degenerate with a 1-dimensional kernel.  Let $\Ker B|_V=\Kee d$, where $d$ is an element in the kernel of the restriction $B|_V$ such that $B(d,c)= 1$. Since the kernel of the invariant form is an ideal, and $d$ can not lie in $\fg^{(1)}$, we see that $V=\fh\oplus\Kee d$, where $\fh$ is an arbitrary subspace of $V$ complementing $\Kee d$ and containing $\fg^{(1)}$, and hence an ideal. Then the Lie (super)algebra $\fg$ is a direct sum of ideals $\fg=\fh\oplus \Kee d\oplus \Kee c$ and the form $B$ is nondegenerate on  $\fh$, as well as on the commutative ideal $\Kee d\oplus \Kee c$. In other words,  this case is \textit{decomposable}.

\sssec{Theorem}\label{Pr5.1} \textit{Let $\fg$ be a Lie algebra over any field  $\Kee$ or a Lie superalgebra over a field $\Kee$ of characteristic $p\ne 2$; let $B$ be a non-degenerate invariant (super)symmetric bilinear form on  $\fg$, and $\fz(\fg)$ the center of $\fg$. If $\fg^{(1)}\cap\fz(\fg)\neq 0$, then $\fg$ is a double extension of a Lie (super)algebra $\fh$.}

\begin{proof} Let $c\ne 0$ be a central element of $\fg$ lying in $\fg^{(1)}$. We have seen already that the space $V:=c^\perp$ contains $\fg^{(1)}$, and hence is an ideal of $\fg$. Since $c\in \fg^{(1)}$, we have $B(c,c)=0$. Non-degeneracy of $B$ implies that, first, $\codim V=1$, and, second, there exists an element $d\in \fg\setminus V$ such that $B(d,c)=1$. Since $c=[x,y]$ for certain $x,y\in \fg$, the invariance of $B$ implies that
$$
-(-1)^{p(x)p(y)}B(y,[x,d])=B([x,y],d)=B(c,d)\ne 0,
$$
i.e., $d$ is not central.

Let $\fh:=V/\Kee c$. The element $c$ belongs to the kernel of the restriction  $B|_V$. Therefore $B$ descends onto the quotient $\fh$, and remains $\fh$-invariant. Denote this restriction by $B_\fh$.  The action of $d$ also descend onto $\fh$, defining a derivation $D$ of $\fh$:
\[
D:=\pr_\fh(d|_V).
\]
Besides, the $\fg$-invariance of $B$ implies that $B_\fh$ is $D$-invariant.

Denote $W:=d^\perp\cap V$. Then $V=\Kee c \oplus W$ (as linear spaces) and the natural projection ${\pr: W\tto \fh}$ is an isomorphism of linear spaces, i.e., we may consider  $W=\pr^{-1}(\fh)$ as \lq\lq embed\-ding of $\fh$ (as a space) into $\fg$".

As a result, we see that $\fg$ is cooked from $\fh$ by means of a central extension and a derivation while the form  $B$ is obtained from the form $B_\fh$ by precisely the rules of constructing double extensions. So it only remains to compute the cocycle  $\sigma$ defining the central extension.

Let $x,y\in W$. Then
$$
\begin{array}{l}
{}[x,y]=\pr^{-1}([\pr x, \pr y])+\sigma(\pr x,\pr y)c\text{~and~}B(d,[x,y])=\\
\sigma(x,y)B(d,c)=\sigma(\pr x,\pr y).
\end{array}$$
Now, we use the invariance condition for the triple $x, d, y$:
$$
\begin{array}{l}
B([x,d],y)+(-1)^{p(x)p(d)}B(d,[x,y])=0\ \Lra \\
B(d,[x,y])=B([d,x], y)=B_\fh(D\pr x,\pr y)\ \Lra\ \
\sigma(\pr x,\pr y)= B_\fh(D\pr x,\pr y).
\end{array}
$$
But this is precisely the statement that $\fg$ is a double extension of $\fh$.
\end{proof}

\ssec{On history} In 1984, the notion of \textit{double extension} of Lie algebras (shorter and more sugges\-tive\-ly called \textit{D-extension} in \cite{BeBou})  was distinguished, see \cite{MR}.  Medina  and Revoy inductively constructed a Lie algebra $\fg$ with a~non-degenerate $\fg$-invariant symmetric bilinear form $B_\fg$ from an algebra $\fh$ of dimension $\dim \fg-2$ with a non-zero center and a non-degenerate $\fh$-invariant symmetric bilinear form $B_\fh$.

At almost the same time there was written a paper \cite{FS}, in which the doubly extended Lie algebras $\fg$ and $\tilde \fg$ were considered up to \textit{isometry}, i.e., an isomorphism $\pi:\fg\rightarrow \tilde \fg$ such that
\[
B_{\tilde \fg}(\pi(f), \pi(g))=B_\fg(f,g)\text{~~for any $f,g\in \fg$}.
\]
It is reasonable to consider \textit{classes of double extensions up to an isometry}, not individual double extensions. This proved useful, e.g., in \cite{BeBou} where several new classes of double extensions explained some results in \cite{BGLL}.

In \cite{BB}, one of the first papers on double extensions, the notion was extended to Lie superalgebras for $p\neq 2$. For a most recent of various generalizations of Theorem \ref{Pr5.1}, see \cite{ABB}.

The paper \cite{BeBou}  gives a review of known examples of  double extensions, but its main, new, and most interesting results are general constructions and examples of double extensions of Lie superalgebras over fields of characteristic $p=2$.


\section{Lie (super)algebras with indecomposable symmetrizable Cartan matrices}

Any Lie (super)algebra $\fg(A)$ with symmetrizable Cartan matrix  $A$ with entries in the ground field $\Kee$ has an invariant (super)symmetric bilinear form; this form is non-degenerate if $A$ is invertible. For a precise definition of Cartan matrix of the Lie algebra in the case where $p>0$ and that of the  Lie superalgebra for any $p$, see \cite{BGL}, where finite-dimensional Lie (super)algebras with indecomposable Cartan matrices are classified  over $\Kee$. Two types of infinite-dimensional Lie (super)algebras with indecomposable $A$ are classified and can be studied closer:

(A) finite-dimensional and of polynomial growth (one stringy and  \lq\lq affine Kac--Moody" (super)algebras,  see \cite{LSS, HS, BGL}) and

(B) one class of exponential growth: \lq\lq almost affine", a.k.a. \lq\lq hyperbolic", Lie algebras and superal\-gebras, see \cite{CCLL}.

We denote the positive elements of the Chevalley basis  (for its definition, see \cite{CCLL}), by $x_i$, the corresponding negative
ones are $y_i$; we set $h_i:=[x_i, y_i]$ for the \textit{generators} $x_i,\ y_i$ only.

\ssec{NIS on $\fg(A)$ with $A$ symmetrizable}\label{Asym} Let $A=DB$, where $D=\diag(\eps_1,\dots, \eps_n)$ and $B=B^T$,  be an $n\times  n$ Cartan matrix, $n<\infty$.   For the proof of existence and uniqueness (up to a scalar factor) of NIS in the non-super case over $\Cee$, see \cite[Th.2.2, p.17]{Kb}, the superization  and generalization for algebraically closed fields of   characteristics $p>2$ is immediate if either $p\neq 2$, or $p=2$ and $\fg_{2\alpha}=0$ for any root $\alpha$; namely, as follows.

Define the (super)symmetric invariant bilinear form $(\cdot, \cdot)$ inductively, starting with Chevalley gene\-rators $x_i, y_i$ of degree $\pm1$ as follows, where indices in $\fg_{\pm1}$ are degrees relative the principal $\Zee$-grading, whereas in $\fg_{\alpha}$ they denote the weights:
\be\label{recipe}
\begin{array}{ll}
\text{On $\fg_0$}:&(h_i, h_j):= B_{i,j}\eps_i\eps_j;\\[1mm]
\text{On $\fg_{\pm1}$}:&(x_i, y_j):= \delta_{i,j}\eps_j;\\[1mm]
\text{On $\fg_{\alpha}\oplus\fg_{-\alpha}$}:&(x, y)= (-1)^{p(x)p(y)}(y, x):= \sum([x, u_i], v_i)\\
&\text{for any $x\in\fg_{\alpha}$ and any representation}\\
&\text{$y=\sum[u_i, v_i]\in \fg_{-\alpha} $, where $0>\deg u_i>\deg y, \  0>\deg v_i>\deg y$};\\[2mm]
\fg_\alpha\perp \fg_\beta&\text{if $\alpha+\beta\neq 0$}.\\
\end{array}
\ee
The form $(\cdot, \cdot)$ is NIS on $\fg(A)$ for any symmetrizable $A$, even if $A=0$. (For example, if $A$ is non-invertible of corank 1,  the form induces NIS on $\fg(A)^{(1)}/\fc$ , where  $\fg(A)^{(1)}:=[\fg(A), \fg(A)]$ and $\fc$ is the center of $\fg(A)$. This, however, is hardly interesting: the point is given a NIS on $\fh$, define a NIS on a double extension of $\fh$.)

In the next subsections we consider \textbf{finite-dimensional} modular Lie (super)algebras $\fg(A)$ with indecomposable Cartan matrix $A$ classified in \cite{BGL}. All such Lie algebras (resp. superal\-geb\-ras) are rigid for $p>3$ (resp. for $p>5$, with the exception of $\fosp_a(4|2)$), see \cite{BGLd}. Below we consider all non-rigid modular  Lie algebras and Lie superalgebras with indecomposable Cartan matrix and their simple relatives, see \cite{BLW} and \cite{BGLd}, except for those deformed with an odd parameter.

\sssbegin{Remark}
It seems that for $p= 2$ Lie superalgebras with root spaces $\fg_{2\alpha}\neq 0$ can not have NIS since even $\fo\fo_{I\Pi}^{(1)}(1|2)$ --- the analog of $\fosp(1|2)$ for $p=2$, see \cite{BGL} --- does not have it. If the form $(\cdot, \cdot)$ is invariant, then nondegeneracy is violated:
\[
\begin{array}{l}
(x,x^2)=([x^2,y],x^2)=(y,[x^2,x^2])=0;\\
(x,x)=([x^2,y],x)=(y,[x^2,x])=0;\\
(x,h)=([x^2,y],h)=(y,[x^2,h])=0;\\
(x, y) = ([x^2, y], y) = (x^2, [y, y]) = 0;\\
(x,y^2)=([x^2,y],y^2)=(x^2,[y,y^2])=0.
\end{array}
\]
 \end{Remark}

\ssec{The deforms of $\fo(5)$ for $p=3$} These deforms constitute  two  non-isomorphic types, see \cite{BLW}:

$1)$ The parametric family $\fbr(2; \eps)$ with Cartan matrix
$\begin{pmatrix} 2 & -1 \\ -2 & 1 - \eps \end{pmatrix}$; it has NIS by the recipe \eqref{recipe}.
Note that $\fbr(2; \eps)$ is simple if $\eps\neq0$, and $\fbr(2;-1)=\fo(5)\simeq\fsp(4)$.

$2)$ An exceptional simple Lie algebra $\mathfrak{L}(2,2)$ discovered by A.~Kostrikin and Kuznetsov. Recall the description of $\mathfrak{L}(2,2)$, see \cite[Prop.~3.2]{BLW}.
The \textit{contact bracket} of two polynomials in divided powers $f,g\in \cO(p,q,t;\un)$ 
is defined to be
\begin{equation}\label{cb}
 \{f,g\}:=\triangle f\cdot\partial_t g - \partial_t f\cdot\triangle g
 +\partial_p f\cdot\partial_q g - \partial_q f\cdot\partial_p g,\text{~~where $\triangle f=2f - p\partial_p f - q\partial_q f$.}
\end{equation}

Let $\alpha$ and $\beta$ be the
simple roots of $\fo(5)$ and $E_\gamma$ be the root vector corresponding
to root $\gamma$. The following is the bracket in  $\mathfrak{L}(2,2)$
\be\label{L-1bracket}
\begin{array}{ccl}
[\cdot,\cdot]&=& \{\cdot,\cdot\}+2\, c, \text{~~where $c=  x_1\otimes dy_3 \wedge dy_4+ 2 \,
   x_3\otimes dy_1 \wedge dy_4 +
   x_4\otimes dy_1 \wedge dy_3$}
\end{array}
\ee
given in terms a basis of $\fbr(2; \eps)$  expressed via generating functions of $\mathfrak{k}(3; \One)$:
\begin{equation}\label{tab}\small
\renewcommand{\arraystretch}{1.4}
\begin{tabular}{|c|l|} \hline
$\deg$&the generator${}_{\text{its weight}}\sim$ its generating function(=the Chevalley basis vector) \\
\hline \hline
$-2$&$E_{-2\alpha -\beta}= \{E_{-\alpha}, E_{-\alpha-\beta}\}\sim 1(=y_4)$\\
\hline
$-1$& $E_{-\alpha}\sim p(=y_2);\quad E_{-\alpha-\beta}= \{E_{-\beta}, E_{-\alpha}\}\sim q(=y_3)$  \\
\hline $0$& $H_{\alpha}\sim - \eps t+ pq(=h_2);\quad H_{\beta}\sim -pq(=h_1); \quad
E_{\beta}\sim p^2(=y_1); \quad
E_{-\beta}\sim -q^2(=x_1)$  \\
\hline $1$& $E_{\alpha}\sim -(1+\eps)pq^2+ \eps
qt(=x_2);\quad E_{\alpha+\beta}= \{E_\beta, E_\alpha\}\sim (1+\eps)p^2q+\eps p t(=x_3)$  \\
\hline $2$& $E_{2\alpha+\beta}= \{E_\alpha,
E_{\alpha+\beta}\}\sim \eps(1+\eps)p^2q^2+ \eps^2t^2(=x_4)$  \\
\hline
\end{tabular}
\end{equation}
\normalsize

\sssec{Claim (NIS on deforms of $\fo(5)$)}\label{cl_o(5)}\textit{
The deform~\eqref{L-1bracket} of $\fo(5)$ preserves NIS with the same Gram matrix as the one determined by recipe~\eqref{recipe} applied to~$\fo(5)$.}

In Claims \ref{wk4}, \ref{wk3} we consider all Lie (super)algebras with indecomposable Cartan matrix (not isomorphic to $\fbr(2; \eps)$) that can be deformed with even parameter, see \cite{BGLd}, but whose deforms have no Cartan matrix. We consider deforms with even parameter as the same space but with different brackets; in this way we can compare the Gram matrices of bilinear forms on the algebra and its deform.

\ssec{Claim (NIS on deforms of $\fg=\fbr(3)$ and $\fwk(4;\alpha)$)}\label{wk4} \textit{For Lie algebras $\fg=\fbr(3)$, and $\fwk(4;\alpha)$ for $\alpha\neq 0, 1$, where $\fwk(4;\alpha)\simeq\fwk(4;\alpha')$ if and only if $\alpha'= \nfrac{1}{\alpha}$, see \cite{BGL}, and Lie superalgebra $\fbrj(2;3)$, all deforms depend on even parameters \textup{(}as proved in \cite{BGLd}\textup{)}. These deforms preserve NIS with the same Gram matrix as the one determined by recipe~\eqref{recipe} applied to~$\fg$, except for the deform of $\fwk(4;\alpha)$ with cocycle $\lambda c_0$ with  parameter $\lambda\in\Kee$ distinct from $0$ and $\alpha$, when the Gram matrix is a different one, of the form
$\begin{pmatrix} \Gamma_{h,h}&0&0\\
0&0&\Gamma_{x,y}\\
0&(\Gamma_{x,y})^T&0\end{pmatrix}$, where for the same numbering of Chevalley generators $h$, $x$ and $y$
as in \cite{BGL} for the Cartan matrix obtained from~ \eqref{CM} at $ \lambda=0$, we have
\be\label{CM}
\Gamma_{h,h}=\begin{pmatrix}
0 & \alpha  & 1 & \lambda  \cr
\alpha  & 0 & 0 & 0 \cr
1 & 0 & 0 & 1 \cr
\lambda  & 0 &1 & 0 \cr
\end{pmatrix}
\ee
and
\be\label{CM1}
\Gamma_{x,y}=\diag\left(1, \ \nfrac{\alpha }{\alpha +\lambda },\  1,\ 1, \ \alpha,\  1,\ 1,\ \alpha,\ 1, \  A_2, \  \alpha, \  A_2 , \  A_2 , \  A_3, \ \alpha A_3\right),
\ee
where $A_2:=\left(\alpha ^2+\lambda  \alpha
+\alpha \right)$ and $A_3:=\left(\alpha
^3+\alpha ^2+\lambda ^2 \alpha +\lambda  \alpha \right)$.}

\textit{For the deform of~$\fwk(4;\alpha)$ with cocycle~$\alpha c_0$, i.e., when the deformation parameter is
equal to~$\alpha$, we only have a \textbf{degenerate} invariant
symmetric bilinear form $B_{\alpha c_0}$ for which $B_{\alpha c_0}(x_2, y_2)\neq 0$ and zero on all other pairs of Chevalley basis vectors.}

\ssec{Claim (NIS on deforms of $\fwk(3;\alpha)$)}\label{wk3}  \textit{For Lie algebra $\fwk(3;\alpha)$ for $\alpha\neq 0,
1$, where $\fwk(3;\alpha)\simeq\fwk(3;\alpha')$ if and only if $\alpha'=\nfrac{a\alpha +b}{c\alpha+d}$ for some $\begin{pmatrix} a&b\\
c&d\end{pmatrix}\in\SL(2;\Zee/2)$, see \cite{BGL}, all deforms depend on even parameters
\textup{(}\cite{BGLd}\textup{)}.
These deforms preserve NIS with the same Gram matrix as the one
determined by recipe~\eqref{recipe} applied to~$\fg$ except for the
deform with cocycle $\lambda c_0$ with parameter $\lambda\in\Kee$ distinct from $0$ and $\alpha$, when
the Gram matrix is a different one, of the form
$B_{\alpha,\lambda}:=\begin{pmatrix}\Delta_{h,h}&0&0\\
0&0&\Delta_{x,y}\\
0&(\Delta_{x,y})^T&0\end{pmatrix}$, where for the same numbering of
Chevalley generators $h$, $x$ and $y$
as in \cite{BGL} for the Cartan matrix obtained from the  matrix~\eqref{CMwk3def} at $ \lambda=0$, we have
\be\label{CMwk3def}
\Delta_{h,h}=\begin{pmatrix}
0 & \alpha +\lambda  & 0 \\
\alpha +\lambda  & 0 & \nfrac{\alpha +\lambda }{\alpha } \\
0 & \nfrac{\alpha +\lambda }{\alpha } & 0
\end{pmatrix}
\ee
and
\[
\Delta_{x,y}=\diag\left(
1,\ \nfrac{A }{\alpha},\ \nfrac{A }{\alpha },\ A,\ \nfrac{A }{\alpha }, \ A,\ \alpha^2 \lambda ^2+\alpha^2+\alpha  \lambda^3+\alpha  \lambda^2+\alpha
\lambda +\alpha +\lambda^3+\lambda
\right),
\]
where $A=\alpha\lambda +\alpha +\lambda^2+\lambda$.}

\textit{The kernel of $B_{\alpha,\lambda}$ is the center of $\fwk^{(1)}(3;\alpha)$ spanned by the vector $c:=h_1 + \alpha h_3$ and the restriction of $B_{\alpha,\lambda}$ to $\fwk^{(1)}(3;\alpha)/\Kee c$ is a NIS.}

\textit{For the deform of~$\fwk^{(1)}(3;\alpha)/\Kee c$ with cocycle~$\alpha c_0$, i.e., when the deformation parameter is
equal to~$\alpha$, we only have a \textbf{degenerate} invariant
symmetric bilinear form $B_{\alpha c_0}$ for which $B_{\alpha c_0}(x_1, y_1)\neq 0$ and zero on all other pairs of Chevalley basis vectors.}

\section{Linear (matrix) and vectorial Lie (super)algebras; their simple relatives}

In representation theory, it is reasonable to consider Lie (super)algebras $\fg$ of vector fields together with a natural topology defined by Weisfeiler filtration determined, in its turn, by a grading vector $r$, see \cite{LSh, BGLLS}. In our search for NIS on $\fg$, topology on the space of $\fg$ is irrelevant, and we consider these Lie (super)algebras as abstract. We recall the list of all known simple finite-dimensional modular Lie (super)algebras for $p\geq 3$, as well as the list of simple Lie (super)algebras of polynomial growth over $\Cee$, see \cite{LSh, GLS}; we indicate which of them have NIS.

We describe each vectorial Lie superalgebra as the result of the (generalized) Cartan prolongation of a pair consisting of a Lie superalgebra and a module over it.    

\ssec{Linear (matrix) Lie (super)algebras}\label{SS:2.3} The \textit{general linear} Lie
superalgebra of all supermatrices of size $\Par$ corresponding to linear operators in the superspace $V=V_{\bar 0}\oplus V_{\bar 1}$ over the ground field $\Kee$ is denoted by
$\fgl(\Par)$, where $\Par=(p_1, \dots, p_{|\Par|})$ is an ordered
collection of parities of the basis vectors of $V$ and $|\Par|=\dim V$;
usually, for the \textit{standard} (simplest) format, $\fgl(\ev,
\dots, \ev, \od, \dots, \od)$ is abbreviated to $\fgl(\dim V_{\bar
0}|\dim V_{\bar 1})$. Any supermatrix from $\fgl(\Par)$ can
be uniquely expressed as the sum of its even and odd parts; in the
standard format this is the following block expression; on non-zero summands the parity is defined:
\[
\mmat A,B,C,D,=\mmat A,0,0,D,+\mmat 0,B,C,0,,\quad
 p\left(\mmat A,0,0,D,\right)=\ev, \; p\left(\mmat 0,B,C,0,\right)=\od.
\]

The \textit{supertrace} is the map $\fgl (\Par)\tto \Kee$,
$(X_{ij})\longmapsto \sum (-1)^{p_{i}}X_{ii}$, where $\Par=(p_{1},
\dots, p_{|\Par|})$. Thus, in the standard format, $\str \begin{pmatrix}A&B\\ C&D\end{pmatrix}=\tr A- \tr D$.

Observe that for Lie superalgebra $\fgl_\cC(p|q)$ over a supercommutative superalgebra $\cC$ we have 
\[
\str X=\tr A- (-1)^{p(X)}\tr D\text{~~for $X=\begin{pmatrix}A&B\\ C&D\end{pmatrix}$,}
\]
so on odd supermatrices with entries in $\cC$ such that $\cC_\od\neq 0$, the supertrace coincides with the trace.

Since $\str [x, y]=0$, the subsuperspace of supertraceless
matrices constitutes the \textit{special linear} Lie subsuperalgebra
$\fsl(\Par)$.

There are, however, at least two super versions of $\fgl(n)$, not
one; for reasons, see \cite[Ch1, Ch.7]{SoS}. The other version --- $\fq(n)$ --- is called the \textit{queer}
Lie superalgebra and is defined as the one that preserves the
complex structure given by an \textit{odd} operator $J$, i.e.,
$\fq(n)$ is the centralizer $C(J)$ of $J$:
\[
\fq(n)=C(J)=\{X\in\fgl(n|n)\mid [X, J]=0 \}, \text{ where }
J^2=-\id.
\]
It is clear that by a~change of basis we can reduce $J$ to the form (shape)
\[
J_{2n}=\antidiag(1_n,-1_n):=\mat {0&1_n\\-1_n&0}
\]
in the standard
format and then $\fq(n)$ takes the form
\begin{equation}\label{q}
\fq(n)=\left \{(A,B):=\mat {A&B\\B&A}, \text{~~where $A, B\in\fgl(n)$}\right\}.
\end{equation}

On $\fq(n)$, the \textit{queertrace} is defined: $\qtr\colon (A,B)\longmapsto \tr B$. Denote by $\fsq(n)$ the Lie superalgebra
of \textit{queertraceless} matrices; set $\fp\fsq(n):=\fsq(n)/\Kee 1_{2n}$.

\sssec{Supermatrices of operators}
To the linear map of superspaces $F: V\tto W$ there corresponds the
dual map $F^*:W^*\tto V^*$ between the dual superspaces. In bases
consisting of the homogeneous vectors $v_{i}\in V$ of parity $p(v_i)$, and $w_{j}\in W$ of parity $p(w_j)$, the formula
$F(v_{j})=\mathop{\sum}_{i}w_{i}X_{ij}$ assigns to $F$ the
supermatrix $X$. In the dual bases, the \textit{supertransposed}
matrix $X^{st}$ corresponds to $F^*$:
\[
(X^{st})_{ij}=(-1)^{(p(v_{i})+p(w_{j}))p(w_{j})}X_{ji}.
\]

\sssec{Supermatrices  of bilinear forms}\label{sssBilMatr}
The supermatrices $X\in\fgl(\Par)$ such that
\[
X^{st}B+(-1)^{p(X)p(B)}BX=0\quad \text{for an homogeneous matrix
$B\in\fgl(\Par)$}
\]
constitute the Lie superalgebra $\faut (B)$ that preserves the
bilinear form $B^f$ on $V$ whose matrix $B=(B_{ij})$ is given by the formula
\be\label{martBil}
B_{ij}=(-1)^{p(B)p(v_i)}B^f(v_{i}, v_{j})\text{~~ for the basis vectors $v_{i}\in V$.}
\ee

In order to identify a bilinear form $B(V, W)$ with an operator, an element of $\Hom(V, W^*)$, the matrix~$B$ of the bilinear form~$B^f$ is defined in \cite[Ch.1]{SoS} by eq.~\eqref{martBil}, not by seemingly natural but inappropriate for such an identification formula
\be\label{wrong}
B_{ij}=B^f(v_i,v_j)\text{~~for the basis vectors $v_i\in V$.}
\ee
Moreover, for the odd forms $B$, the definition \eqref{wrong} contradicts the obvious symmetry of the NIS defined by $\qtr$ on $\fq$. Indeed, the \textit{symmetry} of a homogeneous form $B^f$ means, according to \cite[Ch.1]{SoS}, that $B^f(v, w)=(-1)^{p(v)p(w)}B^f(w,v)$ for any $v \in V$ and $w\in W$, i.e.,
its matrix $B=\mmat R,S,T,U,$ satisfies the condition
\be\label{BilSy}
B^{u}=B,\;\text{ where $B^{u}=
\mmat R^{t},(-1)^{p(B)}T^{t},(-1)^{p(B)}S^{t},-U^{t},$.}
\ee
Similarly, \textit{antisymmetry} of $B$ means that $B^{u}=-B$.
Thus, we see that the \textit{upsetting} of bilinear forms
$u\colon\Bil (V, W)\tto\Bil(W, V)$, which for the \textit{spaces}
and the case where $V=W$ is expressed on matrices in terms of the
transposition, is a~new operation, \textit{not} supertransposition.

Observe
that \textbf{the passage from $V$ to $\Pi (V)$ turns every symmetric
form $B$ on $V$ into an antisymmetric one on $\Pi (V)$}.

Most popular normal forms (shapes) of the even nondegenerate supersymmetric
form are the ones whose supermatrices in the standard format are in the
following normal forms:
\[
\begin{array}{l}
B_{ev}(m|2n)= \diag(1_m, J_{2n}):=\mmat 1_m,0,0,J_{2n}, \text{~or~} \diag(A_m,J_{2n}):=\mmat A_m,0,0,J_{2n},,\\ \text{where
$J_{2n}= \antidiag(1_n,-1_n):=\mmat 0, 1_n,-1_{n},0,$ and $A_m=\antidiag(1,\dots, 1)$}.\\
\end{array}
\]

The usual notation for $\faut (B_{ev}(m|2n))$ is $\fosp(m|2n)$;
sometimes one writes more explicitly, $\fosp^{sy}(m|2n)$. Observe
that the antisymmetric non-degenerate bilinear form is preserved by the
``symplectico-orthogonal" Lie superalgebra, $\fsp\fo (2n|m)$ or,
more prudently, $\fosp^{a}(m|2n)$, which is isomorphic to
$\fosp^{sy}(m|2n)$.

A nondegenerate \textbf{symmetric} odd bilinear form $B_{odd}(n|n)$ can
be reduced to a~normal shape whose matrix in the standard format is
$J_{2n}$, see \eqref{BilSy}, NOT $\Pi_{2n}:=\antidiag(1_n,1_n)$. A normal shape of the \textbf{anti}symmetric odd nondegenerate
form in the standard format is $\Pi_{2n}$. The usual notation for
$\faut (B_{odd}(\Par))$ is $\fpe(\Par)$.
The passage from $V$ to $\Pi (V)$ establishes an isomorphism
$\fpe^{sy}(\Par)\cong\fpe^{a}(\Par)$. These isomorphic Lie
superalgebras are called, as A.~Weil suggested, \textit{periplectic}.

Observe that, though the Lie superalgebras $\fosp^{sy} (m|2n)$ and
$\fosp ^{a} (m|2n)$, as well as $\fpe ^{sy} (n)$ and $\fpe ^{a}
(n)$, are isomorphic, the difference between them is sometimes
crucial, e.g., their Cartan prolongs, see Subsections \ref{2.6},
\ref{SS:2.6.1}, \ref{SS:2.6.3}, are totally different, see
\cite{Sh5}.

The \textit{special periplectic} superalgebra is simple; it is defined to be
\[
\fspe(n)=\{X\in\fpe(n)\mid \str
X=0\}.
\]

Of particular interest to us will also be the Lie superalgebras (here $\Char\Kee>2$)
\begin{equation}
\label{spe}
 \fspe(n)_{a, b}=\fspe(n)\ltimes \Kee(az+bd), \text{
where $z=1_{2n}$, $d=\diag(1_{n}, -1_{n})$, and $a,b\in\Kee$},
\end{equation}
and the nontrivial central extension $\fa\fs$ of $\fspe(4)$ that
we describe after some preparation.

Finally, observe that the term \textit{super}symmetric applied to the bilinear forms in the title of this paper refers to the property $B_{ij}=(-1)^{p(v_i)p(v_j)}B_{ji}$ of the matrix of the bilinear form \eqref{martBil}; according to the definition \eqref{BilSy} this form is \textit{symmetric}.

\ssec{A.~Sergeev's central extension $\fa\fs$ of $\fspe(4)$}\label{ssas}
In 1970's A.~Sergeev proved that over $\Cee$ there is just one
nontrivial central extension of $\fspe(n)$ for $n>2$. It exists only
for $n=4$ and we denote it by $\fas$. (For a generalization of Sergeev's result to analogs of $\fspe(n)$ over fields $\Kee$ of characteristic $p>0$, see \cite{BGLL}.) Let us represent an arbitrary
element $A\in\fas$ as a~pair $A=x+d\cdot z$, where $x\in\fspe(4)$,
$d\in{\Cee}$ and $z$ is the central element. The bracket in $\fas$
is
\begin{equation}
\label{2.1.4}
\left[\mat {a &b\\ c&-a^t} +d\cdot z,
\mat{a' & b' \cr c' & -a'{}^t} +d'\cdot
z\right]= \left[\mat{ a& b \cr c & -a^t},
\mat{ a' & b' \cr c' & -a'{}^t
}\right]+\tr~c\widetilde c'\cdot z,
\end{equation}
where $\ \widetilde {}\ $ is extended via linearity from matrices
$c_{ij}=E_{ij}-E_{ji}$ on which $\widetilde c_{ij}=c_{kl}$ for any
even permutation $(1234)\longmapsto(ijkl)$.

The Lie superalgebra $\fas$ can also be described with the help of
the spinor representation, see \cite{LShs}. For this, consider the Poisson superalgebra $\fpo(0|6)$, the Lie superalgebra whose
superspace is the Grassmann superalgebra $\Lambda(\xi, \eta)$
generated by $\xi_1, \xi_2, \xi_3, \eta _1, \eta _2, \eta_3$ and the
bracket is the Poisson bracket \eqref{2.3.6}.

Recall that $\fh(0|6)=\Span (H_f\mid f\in\Lambda (\xi, \eta))$. Now,
observe that $\fspe(4)$ can be embedded into $\fh(0|6)$, see \cite{ShM}. Indeed,
setting $\deg \xi_i=\deg \eta _i=1$ for all $i$ we introduce
a~$\Zee$-grading on $\Lambda(\xi, \eta)$ which, in turn, induces
a~$\Zee$-grading on $\fh(0|6)$ of the form
$\fh(0|6)=\mathop{\oplus}_{i\geq -1}\fh(0|6)_i$. Since
$\fsl(4)\cong\fo(6)$, we can identify $\fspe(4)_0$ with
$\fh(0|6)_0$.

It is not difficult to see that the elements of degree $-1$ in the
standard gradings of $\fspe(4)$ and $\fh(0|6)$ constitute isomorphic
$\fsl(4)\cong\fo(6)$-modules. It is subject to a~direct verification
that it is possible to embed $\fspe(4)_1$ into $\fh(0|6)_1$.

Sergeev's extension $\fas$ is the result of the restriction of the cocycle that turns $\fh(0|6)$ into
$\fpo(0|6)$ to
$\fspe(4)\subset\fh(0|6)$. The quantization deforms $\fpo(0|6)$ into
$\fgl(\Lambda(\xi))$; the through maps $T_\lambda:
\fas\tto\fpo(0|6)\tto\fgl(\Lambda (\xi))$ are representations of
$\fas$ in the $4|4$-dimensional modules $\spin_\lambda$ isomorphic
to each other for all $\lambda\neq 0$. The explicit form of
$T_\lambda$ is as follows:
\begin{equation}
\label{2.1.4.2} T_\lambda\colon \mat{ a~& b \cr c & -a^t }+d\cdot
z\longmapsto \mat{ a~& b-\lambda \widetilde c \cr c & -a^t
}+\lambda d\cdot 1_{4|4}, 
\end{equation}
where $1_{4|4}$ is the unit matrix and $\widetilde c$ is defined in
the line under eq.~\eqref{2.1.4}. Clearly, $T_\lambda$ is an
irreducible representation for any $\lambda$.

\ssec{Projectivization}\label{SS:2.3.1} If $\fs$ is a~Lie algebra of scalar matrices, and $\fg\subset \fgl (n|n)$ is a~Lie subsuperalgebra containing $\fs$, then the \textit{projective}
Lie superalgebra of type $\fg$ is $\fpg= \fg/\fs$. Examples: $\fpq (n)$, $\fpsq (n)$; $\fpgl (n|n)$, $\fpsl (n|n)$;
whereas $\fpgl (p|q)\cong \fsl (p|q)$ if $p\neq q$.

\ssec{\lq\lq Classical" series of vectorial Lie superalgebra over $\Cee$}
In the table, FD marks the particular cases of finite dimension.
\begin{equation}\label{1.4}\footnotesize
\begin{array}{ll}
\renewcommand{\arraystretch}{1.4}
\begin{tabular}{|l|l|}
\hline $N$&the algebra, conditions for its simplicity\cr \hline

1&$\fvect(m|n)$ for $m\geq 1$  (FD:
$m=0$,\ $n>1$)  \cr \hline

2&$\fsvect(m|n)$ for $m>1$ (FD:
$m=0$,\  $n>2$)  \cr \hline

3&$\fsvect^{(1)}(1|n)$ for $n>1$ \cr \hline

4, FD &$\widetilde{\fsvect}(0|n)$ for $n>2$ \cr \hline

\hline 
5&$\fk(2m+1|n)$ \cr \hline

6&$\fh(2m|n)$ for $m>0$ \cr
\hline

7&$\fh_{\lambda}(2|2)$ for $\lambda\neq -2, -\frac32, -1, \frac12,
0, 1, \infty$ \cr \hline

8, FD&$\fh^{(1)}(0|n)$ for $n>3$ \cr \hline
\end{tabular}
&
\renewcommand{\arraystretch}{1.4}
\begin{tabular}{|l|l|}
\hline $N$&the algebra, conditions for its simplicity\cr \hline

9&$\fm(n)$ \cr
 \hline

10&$\fsm(n)$ for $n>1$ but $n\neq 3$ \cr \hline

11&$\fb_{\lambda}(n)$ for $n>1$ and $\lambda\neq 0, 1,
\infty$ \cr

\hline

12&$\fb_{1}^{(1)}(n)$ for $n>1$ \cr \hline

13&$\fb_{\infty}^{(1)}(n)$ for $n>1$ \cr \hline

14&$\fle(n)$ for $n>1$ \cr
\hline

15&$\fsle^{(1)}(n)$ for $n>2$ \cr
\hline

16&$\widetilde{\fs\fb}_{\mu}(2^{n-1}-1|2^{n-1})$ for $\mu\neq 0$ and
$n>2$\cr \hline
\end{tabular}\\
\end{array}
\end{equation}


We continue explaining the notation used in Table \eqref{1.4} till  Subsection ~\ref{thdefb}.

\textbf{1) General algebras}. Let $x=(u_1, \dots
, u_n, \theta_1, \dots , \theta_m)$, where the $u_i$ are even
indeterminates and the $\theta_j$ are odd ones. Set $\fvect
(n|m):=\fder\; \Cee[x]$; it is called \textit{the general vectorial
Lie superalgebra}.

\textbf{2) Special algebras}. The
\textit{divergence} of the field
$D=\mathop{\sum}_if_i\pder{u_{i}} + \mathop{\sum}_j
g_j\pder{\theta_{j}}$ is the function (in our case: a~polynomial, or
a~series)
\begin{equation}
\label{2.2.4} \Div D=\mathop{\sum}\limits_i\pderf{f_{i}}{u_{i}}+
\mathop{\sum}\limits_j (-1)^{p(g_{j})}
\pderf{g_{i}}{\theta_{j}}.
\end{equation}

$\bullet$ The Lie superalgebra $\fsvect (n|m):=\{D \in \fvect
(n|m)\mid \Div D=0\}$ is called the \textit{special} (or
\textit{divergence-free}) \textit{vectorial superalgebra}.
Equivalently,
\[
\fsvect(n|m)=\{ D\in \fvect (n|m)\mid L_D\vvol _x=0\},
\]
where $\vvol_x$ is the
volume form with constant coefficients in coordinates $x$ and $L_D$
the Lie derivative along the vector field $D$.

$\bullet$ The Lie superalgebra $\fsvect^{(1)}(1|m):=[\fsvect(1|m),
\fsvect(1|m)]$ is called the \textit{traceless special vectorial
superalgebra}.

$\bullet$ The deform of $\fsvect(0|m)$ is the Lie superalgebra
\[
\fsvect_{\lambda}(0|m):=\{D \in \fvect (0|m)\mid \Div
(1+\lambda\theta_1\cdot \dots \cdot \theta_m)D=0\},
\]
where $p(\lambda)\equiv m\pmod 2$. It is called the \textit{deformed
special} (or \textit{divergence-free}) \textit{vectorial
superalgebra}. Clearly, $\fsvect_{\lambda}(0|m)\cong
\fsvect_{\mu}(0|m)$ for $\lambda\mu\neq 0$. So we briefly denote
these deforms by $\widetilde{\fsvect}(0|m)$.
Observe that for $m$ odd the parameter of deformation, $\lambda$, is
odd.

\textbf{3) The algebras preserving Pfaff equations and certain
differential 1- and 2-forms}.

$\bullet$ Set $u=(t, p_1, \dots , p_n, q_1, \dots , q_n)$; let
\begin{equation}
\label{2.2.5} \widetilde \alpha_1 := dt +\mathop{\sum}\limits_{1\leq
i\leq n}(p_idq_i - q_idp_i)\ + \mathop{\sum}\limits_{1\leq j\leq
m}\theta_jd\theta_j\quad\text{and}\quad \widetilde
\omega_0:=d\widetilde
\alpha_1\ . 
\end{equation}
The form $\widetilde \alpha_1$ is called \textit{contact}, the form
$\widetilde \omega_0$ is called \textit{symplectic}. Sometimes
it is more convenient to set
$\Theta=\begin{cases}(\xi, \eta)& \text{ if }\ m=2k\\(\xi, \eta,
\theta)&\text{ if }\ m=2k+1,\end{cases}$ where
\begin{equation}
\label{zam}
\xi_j=\frac{1}{\sqrt{2}}(\theta_{j}-\sqrt{-1}\theta_{k+j}),\quad
\eta_j=\frac{1}{ \sqrt{2}}(\theta_{j}+\sqrt{-1}\theta_{k+j})\; \text{
for}\; j\leq k= \left[\frac{m}{2}\right];\;
\; \theta =\theta_{2k+1}
\end{equation}
and in place of $\widetilde \omega_0$ or $\widetilde \alpha_1$ take
$\alpha_1$ and $\omega_0=d\alpha_1$, respectively, where
\begin{equation}
\label{2.2.6} \alpha_1:=dt+\mathop{\sum}\limits_{1\leq i\leq
n}(p_idq_i-q_idp_i)+ \mathop{\sum}\limits_{1\leq j\leq
k}(\xi_jd\eta_j+\eta_jd\xi_j)
\begin{cases}&
\text{ if }\ m=2k\\
+\theta d\theta&\text{ if }\ m=2k+1.\end{cases} 
\end{equation}

The \textit{contact} Lie superalgebra is the one that preserves the \textit{Pfaff equation}
\[
\text{$\alpha_1(X)=0$  for $X\in \fvect(2n+1|m)$}
\]
or, equivalently, preserves the distribution singled out by the form $\alpha_1$, i.e., the superalgebra
\begin{equation}
\label{k} \fk (2n+1|m)=\{ D\in \fvect (2n+1|m)\mid
L_D\alpha_1=f_D\alpha_1\text{ for some }f_D\in \Cee [t, p, q,
\theta]\}.
\end{equation}
The Lie superalgebra
\begin{equation}
\label{po}
\begin{array}{c}
\fpo (2n|m)=\{ D\in \fk (2n+1|m)\mid L_D\alpha_1=0\}\end{array}
\end{equation}
is called the \textit{Poisson} superalgebra.

The above ``symmetric" expression of $\alpha_1$ is popular among
algebraists; due to its symmetry it is convenient in computations.
In mechanics and differential geometry  (in pre-super era), and in characteristic $2$, the following
expression of the form $\alpha_1$  is used:
\begin{equation}
\label{alphachar2} \alpha_{1(2)}:=dt-\mathop{\sum}\limits_{1\leq
i\leq n}p_idq_i + \mathop{\sum}\limits_{1\leq j\leq k}\xi_jd\eta_j
\begin{cases}&
\text{ if }\ m=2k\\
+\theta d\theta&\text{ if }\ m=2k+1.\end{cases} 
\end{equation}

$\bullet$ Similarly, set $u=q=(q_1, \dots , q_n)$,
let $\theta=(\xi_1, \dots , \xi_n; \tau)$ be odd. Set (expressions in parentheses are for characteristic $2$)
\begin{equation}
\label{alp0}
\begin{array}{c}
\alpha_0:=d\tau+\mathop{\sum}\limits_i(\xi_idq_i+q_id\xi_i), \ \ (\text{or $\alpha_{0(2)}:=d\tau+\mathop{\sum}\limits_i\xi_idq_i$});
\qquad\omega_1:=d\alpha_0 \ \ (\text{resp.  $d\alpha_{0(2)}$})
\end{array}
\end{equation}
and call these forms, as A.~Weil advised, the \textit{pericontact}
and \textit{periplectic}, respectively.

The \textit{pericontact} Lie superalgebra is the one that  preserves the Pfaff equation
\[
\text{$\alpha_0(X)=0$ for $X\in \fvect(n|n+1)$}
\]
or, equivalently, preserves the distribution singled out by the form $\alpha_0$:
\begin{equation}
\label{m} \fm (n)=\{ D\in \fvect (n|n+1)\mid L_D\alpha_0=f_D\cdot
\alpha_0\text{ for some }\; f_D\in \Cee [q, \xi, \tau]\}.
\end{equation}

The Lie superalgebra
\begin{equation}
\label{but}
\begin{array}{c}
\fb (n)=\{ D\in \fm (n)\mid L_D\alpha_0=0\}\end{array}
\end{equation}
is called the \textit{Buttin} superalgebra in honor of C.~Buttin who
was the first to show that the \textit{Schouten bracket} of the functions that generate $\fb (n)$, see \eqref{2.3.7}, satisfies the super
Jacobi identity.

The following are respective \textit{divergence-free} (or \textit{special}) Lie subsuperalgebras
\begin{equation}
\label{sm1}
\begin{array}{c}
\fsm (n)=\{ D\in \fm (n)\mid \Div\ D=0\}, \ \ \
 \fs\fb (n)=\{ D\in
\fb (n)\mid \Div\ D=0\}.
\end{array}
\end{equation}

\ssec{Generating functions}\label{SS:2.4.1}{}~{}

$\bullet$ \underline{Odd form $\alpha_1$}. For any $f\in\Cee [t, p,
q, \theta]$, set:
\begin{equation}
\label{2.3.1} K_f:=(2-E)(f)\pder{t}-H_f + \pderf{f}{t} E,
\end{equation}
where $E:=\mathop{\sum}_i y_i \pder{y_{i}}$ (here the $y_{i}$
are all the coordinates except $t$) is the \textit{Euler operator},
and $H_f$ is the hamiltonian field with Hamiltonian $f$ that
preserves $d\widetilde \alpha_1$:
\begin{equation}
\label{2.3.2} H_f:=\mathop{\sum}\limits_{i\leq n}\left(\pderf{f}{p_i}
\pder{q_i}-\pderf{f}{q_i} \pder{p_i}\right )
-(-1)^{p(f)}\left(\mathop{\sum}\limits_{j\leq m}\pderf{ f}{\theta_j}
\pder{\theta_j}\right ) . 
\end{equation}
The choice of the form $\alpha_1$ instead of $\widetilde \alpha_1$
only affects the shape of $H_f$ that we give for $m=2k+1$:
\begin{equation}
\label{2.3.2'} H_f:=\mathop{\sum}\limits_{i\leq n}\left
(\pderf{f}{p_i} \pder{q_i}-\pderf{f}{q_i} \pder{p_i}\right)
-(-1)^{p(f)}\left(\mathop{\sum}\limits_{j\leq
k}\left(\pderf{f}{\xi_j} \pder{\eta_j}+ \pderf{f}{\eta_j}
\pder{\xi_j}\right)+ \pderf{f}{\theta} \pder{\theta}\right).
\end{equation}

 $\bullet$ \underline{Even form $\alpha_0$}. For any $f\in\Cee [q,
\xi, \tau]$, set:
\begin{equation}
\label{2.3.3} 
M_f:=(2-E)(f)\pder{\tau}- Le_f -(-1)^{p(f)}
\pderf{f}{\tau} E, 
\end{equation}
where $E:=\mathop{\sum}_iy_i \pder{y_i}$ (here the $y_i$ are
all the coordinates except $\tau$), and
\begin{equation}
\label{2.3.4} \Le_f:=\mathop{\sum}\limits_{i\leq n}\left(
\pderf{f}{q_i}\ \pder{\xi_i}+(-1)^{p(f)} \pderf{f}{\xi_i}\
\pder{q_i}\right).
\end{equation}
Since
\begin{equation}
\label{2.3.5}
\renewcommand{\arraystretch}{1.4}
\begin{array}{l}
  L_{K_f}(\alpha_1)=2 \pderf{f}{t}\alpha_1=K_1(f)\alpha_1, \\
L_{M_f}(\alpha_0)=-(-1)^{p(f)}2 \pderf{
f}{\tau}\alpha_0=-(-1)^{p(f)}M_1(f)\alpha_0,
\end{array}
\end{equation}
it follows that $K_f\in \fk (2n+1|m)$ and $M_f\in \fm (n)$.
Observe that
\[
p(\Le_f)=p(M_f)=p(f)+\od.
\]

$\bullet$ To the (super)commutators $[K_f, K_g]$ or $[M_f, M_g]$
there correspond \textit{contact brackets} of the generating functions:
\[
\renewcommand{\arraystretch}{1.4}
\begin{array}{l}
{}[K_f, K_g]=K_{\{f, \; g\}_{k.b.}};\qquad
 {}[M_f, M_g]=M_{\{f, \;
g\}_{m.b.}}.\end{array}
\]
The explicit formulas for the contact brackets are as follows. Let
us first define the brackets on functions that do not depend on $t$
(resp. $\tau$).

The \textit{Poisson bracket} $\{\cdot , \cdot\}_{P.b.}$ (in the
realization with the form $\widetilde{\omega}_0$) is given by the
formula
\begin{equation}
\label{2.3.6}
\renewcommand{\arraystretch}{1.4}
\begin{array}{ll}
\{f, g\}_{P.b.}&:=\mathop{\sum}\limits_{i\leq n}\
\bigg(\displaystyle\pderf{f}{p_i}\
\displaystyle\pderf{g}{q_i}-\displaystyle\pderf{f}{q_i}
\displaystyle\pderf{g}{p_i}\bigg)-\\
&(-1)^{p(f)}\mathop{\sum}\limits_{j\leq m}\
\displaystyle\pderf{f}{\theta_j}\
\displaystyle\pderf{g}{\theta_j}\text{ for any }f, g\in
\Cee [p, q, \theta]\end{array}
\end{equation}
and in the realization with the form $\omega_0$ for $m=2k+1$ it is
given by the formula
\begin{equation}
\label{pb}
\renewcommand{\arraystretch}{1.4}
\begin{array}{ll}
  \{f, g\}_{P.b.}&:=\mathop{\sum}\limits_{i\leq n}\ \bigg(\displaystyle\pderf{f}{p_i}\
\pderf{g}{q_i}-\ \pderf{f}{q_i}\
\pderf{g}{p_i}\bigg)-\\
&(-1)^{p(f)}\bigg(\mathop{\sum}\limits_{j\leq m}\left(
\displaystyle\pderf{f}{\xi_j}\ \pderf{ g}{\eta_j}+\pderf{f}{\eta_j}\
\pderf{ g}{\xi_j}\right)+\displaystyle\pderf{f}{\theta}\ \pderf{
g}{\theta}\bigg)\text{ for any }f, g\in \Cee [p, q, \xi, \eta,
\theta].
\end{array}
\end{equation}

The \textit{Buttin bracket}, or \textit{Schouten bracket}, a.k.a. \textit{antibracket}, $\{\cdot , \cdot\}_{B.b.}$   is given by the formula
\begin{equation}
\label{2.3.7} \{ f, g\}_{B.b.}:=\mathop{\sum}\limits_{i\leq n}\
\bigg(\pderf{f}{q_i}\ \pderf{g}{\xi_i}+(-1)^{p(f)}\
\pderf{f}{\xi_i}\ \pderf{g}{q_i}\bigg)\text{ for any }f, g\in \Cee
[q, \xi].
\end{equation}

 In terms of the Poisson and Buttin brackets,
respectively, the contact brackets are
\begin{equation}
\label{2.3.8} \{ f, g\}_{k.b.}=(2-E) (f)\pderf{g}{t}-\pderf{f}
{t}(2-E) (g)-\{ f, g\}_{P.b.}
\end{equation}
and
\begin{equation}
\label{2.3.9} \{ f, g\}_{m.b.}=(2-E)
(f)\pderf{g}{\tau}+(-1)^{p(f)} \pderf{f}{\tau}(2-E) (g)-\{ f,
g\}_{B.b.}. 
\end{equation}

The Lie superalgebra of \textit{Hamiltonian
fields}  (or \textit{Hamiltonian
superalgebra}) and its special subalgebra (defined only if $n=0$)
are
\begin{equation}
\label{h}
\begin{array}{l}
\fh (2n|m)=\Span(D\in \fvect (2n|m)\mid L_D\omega_0=0),\\
\fh^{(1)} (0|m)=\Span(H_f\in \fh (0|m)\mid \int f\vvol_{\theta}=0).
\end{array}
\end{equation}
The ``odd'' analogues of the Lie superalgebra of Hamiltonian
fields are the Lie superalgebra of vector fields $\Le_{f}$
introduced in \cite{L1} and its special subalgebra:
\begin{equation}
\label{len}
\begin{array}{l}
\fle (n)=\Span(D\in \fvect (n|n)\mid L_D\omega_1=0),\\
\fsle (n)=\Span(D\in \fle (n)\mid \Div D=0).
\end{array}
\end{equation}

It is not difficult to prove the following isomorphisms (as superspaces):
\begin{equation}
\label{iso}
\begin{array}{rcl}
\fk (2n+1|m)&\cong&\Span(K_f\mid f\in \Cee[t, p, q, \xi]),\\
 \fh
(2n|m)&\cong&\Span(H_f\mid f\in \Cee [p, q, \xi]);\\
\fm (n)&\cong&\Span(M_f\mid f\in \Cee [\tau, q, \xi]),\\
\fle
(n)&\cong&\Span(\Le_f\mid f\in \Cee [q, \xi]).\\
\end{array}
\end{equation}
We see that the commutants (first derived algebras) in the purely odd case can be described as
\[
\fpo^{(1)} (0|m)=\Span(K_f\in \fpo (0|m)\mid \int f\vvol_\theta=0);\quad
\fh^{(1)} (0|m)=\fpo^{(1)} (0|m)/\Cee\cdot K_1.
\]

\ssec{The Cartan prolongs}\label{2.6}
We will repeatedly use the Cartan prolongation. So let us recall the
definition and generalize it somewhat. Let $\fg$ be a~Lie algebra,
$V$ a~$\fg$-module, $S^i$ the operator of the $i$th symmetric power.
Set $\fg_{-1} = V$ and $\fg_0 = \fg$.

Recall that, for any (finite-dimensional) vector space $V$, we have
\[
\Hom(V, \Hom(V,\ldots, \Hom(V,V)\ldots)) \simeq L^{k}(V, V, \ldots,
V; V),
\]
where $L^{k}$ is the space of $k$-linear maps and we have
$(k+1)$-many $V$'s on both sides. Now, we recursively define, for
any $k > 0$:
\begin{equation}
\label{prol}
\begin{array}{ll}
 \fg_k = \{X\in \Hom(\fg_{-1}, \fg_{k-1})\mid& X(v_1)(v_2, v_3, ..., v_{k+1})
 =
X(v_2)(v_1, v_3, ..., v_{k+1})\\
&\text{ where $v_1, \dots, v_{k+1}\in \fg_{-1}\}$}. \end{array}
\end{equation}
The space $\fg_k $ is said to be the $k$th \textit{Cartan prolong}
(the result of the Cartan \textit{prolongation}) of the pair
$(\fg_{-1}, \fg_0)$.

Equivalently, let
\begin{equation}
\label{2.5.1}
\begin{array}{l}
i: S^{k+1}(\fg_{-1})^*\otimes \fg_{-1}\tto S^{k}(\fg_{-1})^*\otimes
\fg_{-1}^*\otimes\fg_{-1},\\
j:
S^{k}(\fg_{-1})^*\otimes \fg_{0}\tto
S^{k}(\fg_{-1})^*\otimes \fg_{-1}^*\otimes\fg_{-1}
\end{array}\end{equation}
be the natural maps. Then $\fg_k=i(S^{k+1}(\fg_{-1})^*\otimes
\fg_{-1})\cap j(S^{k}(\fg_{-1})^*\otimes \fg_{0})$.

The \textit{Cartan prolong} of the pair $(V, \fg)$ is $(\fg_{-1},
\fg_{0})_* = \mathop{\oplus}_{k\geq -1} \fg_k$.

If the $\fg_0$-module $\fg_{-1}$ is \textit{faithful},
there exists an injective linear map $\varphi: (\fg_{-1}, \fg_{0})_*\tto \fvect (n)$ such that
\begin{multline}
\varphi((\fg_{-1}, \fg_{0})_*)\subset \fvect (n) = \fder~\Cee[x_1, ... ,
x_n],\; \text{ where }\; n = \dim~\fg_{-1}\; \text{ and }\\ \fg_i
= \{D\in \fvect(n)\mid \deg D=i, [D, X]\in\fg_{i-1}\text{ for any
} X\in\fg_{-1}\}.
\end{multline}
It is subject to a direct verification that the Lie algebra structure
on $\fvect (n)$ induces a~Lie algebra structure on $\varphi((\fg_{-1},
\fg_{0})_*)$. In what follows we do not indicate $\varphi$; the space $(\fg_{-1}, \fg_{0})_*$ has a~natural Lie algebra
structure even if the $\fg_0$-module $\fg_{-1}$ is not faithful.

\ssec{A generalization of the Cartan prolong}\label{SS:2.6.1}
Let
$\fg_-=\mathop{\oplus}_{-d\leq i\leq -1}\fg_i$ be a~
$\Zee$-gra\-ded Lie algebra and $\fg_0\subset \fder_0\fg$ a Lie
subalgebra in the algebra of the $\Zee$-grading-preserving derivations.
Let
\begin{equation}
\label{2.5.1'}
\begin{array}{l} i: S^{k+1}(\fg_{-})^*\otimes \fg_{-}\tto
S^{k}(\fg_{-})^*\otimes \fg_{-}^*\otimes\fg_{-},\\ 
 j: S^{k}(\fg_{-})^*\otimes \fg_{0}\tto
S^{k}(\fg_{-})^*\otimes \fg_{-}^*\otimes\fg_{-}. 
\end{array}\end{equation}
be natural analogs of maps \eqref{2.5.1}. For $k>0$, define the
$k$th prolong of the pair $(\fg_-, \fg_0)$ to be:
\begin{equation}
\label{genprol} \fg_k = \left (j(S^{\bcdot}(\fg_-)^*\otimes
\fg_0)\cap i(S^{\bcdot}(\fg_-)^*\otimes \fg_-)\right )_k,
\end{equation}
where the subscript $k$ in the right hand side singles out the
component of degree $k$.

Set $(\fg_-, \fg_0)_*=\mathop{\oplus}_{k\geq -d} \fg_k$;
then, as is easy to verify, $(\fg_-, \fg_0)_*$ is a~Lie algebra.

Superization of (generalized) Cartan prolongation procedure is immediate.

\ssec{More notation}\label{moreNot} The tautological representation of a~matrix Lie superalgebra $\fg$
or its subalgebra $\fg\subset \fgl(V)$ in $V$ over the ground field $\Kee$, and sometimes \textbf{the module $V$ itself} are denoted by $\id$ or, for clarity,
$\id_{\fg}$. The context prevents confusion of these notations with
that of the identity (scalar) operator $\id_{V}$ on the space $V$,
as in the next paragraph:

For $\fg=\mathop{\oplus}_{i\in\Zee}\fg_{i}$, the trivial
representation of $\fg_0$ is denoted by $\Kee$ (if $\fg_0$ is
simple) whereas $\Kee[k]$ denotes a~representation of $\fg_0$
trivial on the semisimple part of $\fg_0$ and such that $k$ is the value of the
central element $z$ from $\fg_{0}$, where $z$ is chosen so that
$\ad_z|_{\fg_i}=i\cdot \id_{\fg_i}$.

Hereafter, $\fc\fg:=\fg\oplus\Kee z$, the trivial central extension of $\fg$.

\ssec{Vectorial Lie superalgebras as Cartan
prolongs}\label{SS:2.6.2} Superizations of the constructions described in
sec.~\ref{2.6} are straightforward: via Sign Rule. For its application to $p>0$, see Subsetion~\ref{SSlopsup=3}.  We thus get infinite-dimensional Lie superalgebras (superscript ${}^a$ indicates the incarnation of the Lie superalgebra $\fosp$, $\fpe$ or $\fspe$ preserving the anti-symmetric bilinear form):
\begin{equation}
\label{carprol}
\begin{array}{ll}
\begin{array}{l}
\fvect(m|n)=(\id, \fgl(m|n))_*;\\
\fsvect(m|n)=(\id, \fsl(m|n))_*; \\
\fh(2m|n)=(\id, \fosp^{a}(m|2n))_*; \\
\end{array} &\begin{array}{l}
\fle(n)=(\id, \fpe^{a}(n))_*;\\

\fs\fle(n)=(\id, \fspe^{a}(n))_*.\\
\end{array}\end{array}
\end{equation}

The contact Lie superalgebras and exceptional ones are Cartan prolongs $\fg$ with nilpotent (and noncommutative) negative part $\fg_-$;  let us describe them.

$\bullet$ Define the Lie superalgebra $\fhei(2n|m)=W\oplus \Cee z$, where $\dim W=2n|m$ and $W$ is endowed with
a~nondegenerate antisupersymmetric bilinear form $B$, with the following relations:
\begin{equation}
\label{hei} \text{$z$ spans the center; $[v,
w]=B(v, w)\cdot z$ for any $v, w\in W$.}
\end{equation}
Clearly, we have
\begin{equation}
\label{kont} \fk(2n+1|m)=(\fhei(2n|m), \fc\fosp^{a}(m|2n))_*\ .
\end{equation}

$\bullet$ The ``odd'' analog of $\fk$ is associated with the
following ``odd'' analog of $\fhei(2n|m)$. Denote by $\fb\fa(n)=W\oplus \Cee\cdot z$ the
\textit{antibracket} Lie superalgebra ($\fb\fa$ is Anti-Bracket read
backwards), where $W$ is an
$n|n$-dimensional superspace endowed with a~nondegenerate
antisupersymmetric odd bilinear form $B$; the bracket in $\fb\fa(n)$ is
given by the following relations:
\begin{equation}
\label{2.5.4} \text{$z$ is odd and spans the center; $[v,
w]=B(v, w)\cdot z$ for any $v, w\in W$.}
\end{equation}

Clearly,
\begin{equation}
\label{mprol} \fm(n)=(\fb\fa(n), \fc\fpe^{a}(n))_*\ .
\end{equation}


\ssec{A partial Cartan prolong: prolongation of
a positive part}\label{SS:2.6.3}  Let $\fh_1\subset \fg_1$ be
a~$\fg_0$-submodule such that $[\fg_{-1}, \fh_1]=\fg_0$. If such
$\fh_1$ exists (usually, $[\fg_{-1}, \fh_1]\subset\fg_0$), define
the 2nd prolong of $(\mathop{\oplus}_{i\leq 0}\fg_i,
\fh_1)$ to be
\begin{equation}
\label{partprol} \fh_{2}=\{D\in\fg_{2}\mid [D, \fg_{-1}]\subset
\fh_1\}.
\end{equation}
The terms $\fh_{i}$, where $i>2$, are similarly defined. Set
$\fh_i=\fg_i$ for $i\leq 0$ and $\fh_*=\oplus\fh_i$.

\textit{Examples}: $\fvect(1|n; n)$ is a~subalgebra of $\fk(1|2n;
n)$. The former is obtained as the Cartan prolong of the same
nonpositive part as $\fk(1|2n; n)$ and a~submodule of $\fk(1|2n;
n)_1$. The simple exceptional superalgebra $\fk\fas$ introduced in
\cite{Sh5, Sh14}, see table~\eqref{nonpos},  is
another example.

\ssec{Traces and divergencies on vectorial Lie
superalgebras}\label{traceDiv} On any Lie algebra $\fg$ over a field
$\Kee$, a \textit{trace} is any linear map $\tr: \fg\tto \Kee$ such
that
\begin{equation}\label{deftr}
\tr (\fg^{(1)})=0.
\end{equation}

Now, let $\fg$ be a $\Zee$-graded vectorial Lie algebra with
$\fg_{-}:=\mathop{\oplus}\limits_{i<0}\fg_i$ generated by
$\fg_{-1}$, and let $\tr$ be a trace on $\fg_0$. Recall that any
$\Zee$-grading of a given vectorial Lie algebra is given by degrees
of the indeterminates, so the space of functions is also
$\Zee$-graded. Let $\cF$ be the superalgebra of ``functions" (divided powers in
indeterminates $u$ on the $m|n$-dimensional superspace, if $p>0$). The \textit{divergence} $\Div\colon\fg\tto\cF$ is a degree-preserving
$\ad_{\fg_{-1}}$-invariant prolongation of the trace satisfying the following conditions, so that $\Div\in Z^1(\fg; \cF)$, i.e., is a cocycle: 
\[
\begin{array}{l}
X_i(\Div D)=\Div([X_i,D])\text{~~ for all elements
$X_i$ that span $\fg_{-1}$};\\
 \Div|_{\fg_0}=\tr;\\
 \Div|_{\fg_{-}}=0.\end{array}
 \]
\sssec{Divergence-free subalgebras}\label{SS:2.4.3} $\bullet$ In $\fvect(a|b)$, this is $\fsvect(a|b)$.

$\bullet$ In $\fk(2n+1|m)$. Since, as is easy to
calculate,
\begin{equation}
\label{div}
 \Div K_f =(2n+2-m)K_1(f),
\end{equation}
it follows that the divergence-free subalgebra of the contact Lie
superalgebra either coincides with the whole algebra (for $m=2n+2$),
or is isomorphic to the Poisson superalgebra $\fpo(2n|m)$. 

$\bullet$ In $\fm(n)$, the
situation is more interesting: in the standard grading of  $\fm(n)$, the codimension of $(\fm(n)_0)^{(1)}$ in $\fm(n)_0$ is equal to 2; hence, there are two linearly independent traces and two cohomologically inequivalent divergences corresponding to these traces. Let $\Div$ be the restriction of the divergence from  $\fvect(n|n+1)$ to $\fm(n)$: 
\begin{equation}
\label{2.4.1} \Div M_f =(-1)^{p(f)}2\left ((1-E)\pderf{f}{\tau} -
\mathop{\sum}\limits_{i\leq n}\frac{\partial^2 f}{\partial q_i
\partial\xi_i}\right ), 
\end{equation}
The divergence-free (relative to $\Div$) subalgebra of the pericontact
superalgebra is
\begin{equation}
\label{sm2} \fsm (n) = \Span\left (M_f \in \fm (n)\mid
(1-E)\pderf{f}{\tau} =\mathop{\sum}\limits_{i\leq n}\frac{\partial^2
f}{\partial q_i
\partial\xi_i}\right ).
\end{equation}
In particular,
\begin{equation}
\label{2.4.2}
 \Div \Le_f = (-1)^{p(f)}2\mathop{\sum}\limits_{i\leq
n}\frac{\partial^2 f}{\partial q_i \partial\xi_i}.
\end{equation}
Set
\[
\fsle (n):=\Span\left(\Le_f\in \fle(n)\mid\Delta(f)=0, \text{~~where~~} \Delta=\mathop{\sum}\limits_{i\leq n}\frac{\partial^2
}{\partial q_i
\partial\xi_i}\right).
\]

The other divergence can be selected to be $\Pty\circ\partial_\tau$, where $\Pty(X)=(-1)^{p(X)}X$. The subalgebra of $\fm(n)$ corresponding to a linear combination of the two divergences is described in Subsection~\ref{SS:2.8}.

\sssec{Traceless subalgebras}\label{SS:2.4.3a} In terms of super $K_0$-functor, the superdimension $a|b$ is  an element $a+b\eps$, where $\eps^2=1$. Lie superalgebras $\fsle (n)$, $\fs\fb (n)$ and $\fsvect (1|n)$
have traceless ideals $\fsle ^{(1)}(n)$, $\fs\fb ^{(1)}(n)$ and $\fsvect
^{(1)}(n)$ defined from the exact sequences
\begin{equation}
\label{seqv}
\begin{array}{c}
0\tto \fsle ^{(1)}(n)\tto \fsle
(n)\tto \Cee\cdot Le_{\xi_1\dots\xi_n} \tto 0, \ \codim\fsle ^{(1)}(n)= \eps^{n+1};\\
0\tto \fs\fb ^{(1)}(n)\tto \fs\fb
(n)\tto \Cee\cdot M_{\xi_1\dots\xi_n} \tto 0, \ \codim\fs\fb ^{(1)}(n)= \eps^{n+1};\\
\displaystyle 0\tto \fsvect ^{(1)}(n)\tto
\fsvect (1|n)\tto \Cee
\cdot\xi_1\dots\xi_n\partial_{t}\tto 0, \ \codim\fsvect ^{(1)}(n)= \eps^n.\end{array}
\end{equation}

For more examples of traceless subalgebras, see \eqref{2.7.5}.

\ssec{The exceptional simple vectorial Lie superalgebras}
The five exceptional simple vectorial Lie superalgebras are
given below
as Cartan prolongs $(\fg_{-1},
\fg_{0})_{*}$ or generalized  Cartan prolongs $(\fg_{-},
\fg_{0})_{*}$.

For \underline{depth $\leq 2$}, for
$\fg_{-}=\mathop{\oplus}_{-2\leq i\leq -1}\fg_{i}$, we write
$(\fg_{-2}, \fg_{-1}, \fg_{0})_{*}$ instead of $(\fg_{-},
\fg_{0})_{*}$. 
The
corresponding terms $\fg_{i}$ for $i\leq 0$ are given in
~\eqref{nonpos} and \eqref{depth3}; for notation, see Subsection~\ref{moreNot}.

The non-positive components of exceptional simple vectorial Lie superalgebras over
$\Cee$:
\be\label{nonpos}\tiny
\renewcommand{\arraystretch}{1.3}
\begin{tabular}{|c|c|c|c|}
\hline $\fg$&$\fg_{-2}$&$\fg_{-1}$&$\fg_0$\cr
\hline \hline

$\fv\fle(3|6)$&$\id_{\fsl(3)}\boxtimes\Cee_{\fsl(2)}\boxtimes
\Cee[-2]$&$\id^*_{\fsl(3)}\boxtimes \id_{\fsl(2)}\boxtimes
\Cee[-1]$&$\fsl(3)\oplus\fsl(2)\oplus\Cee z$\cr
\hline

$\fv\fas(4|4)$&$-$&$\spin_\lambda$ for any $\lambda\in \Cee^\times$&$\fas$\cr \hline

$\fk\fas(1|6)$&$\Cee[-2]\boxtimes\Cee_{\fo(6)}$&$\Cee\boxtimes\Pi(\id_{\fo(6)})$& $\fc\fo(6)$\cr
\hline

$\fm\fb(3|8)$&$\id_{\fsl(3)}\boxtimes\Cee_{\fsl(2)}\boxtimes
\Cee[-2]$&$\Pi(\id^*_{\fsl(3)}\boxtimes \id_{\fsl(2)}\boxtimes
\Cee[-1])$&$\fsl(3)\oplus\fsl(2)\oplus\Cee z$\cr

\hline

$\fk\fle(5|10)$&$\id$&$\Pi(\Lambda^2(\id^*))$&$\fsl(5)$\cr
\hline
\end{tabular}
\ee
None of the above $\Zee$-graded vectorial Lie
superalgebras \eqref{1.4} and \eqref{nonpos} is of depth $>3$ and only one
is of depth 3:
\begin{equation}\label{depth3}
\fm\fb(3|8)_{-3}=\Pi(\Cee_{\fsl(3)}\boxtimes\id_{\fsl(2)}\boxtimes\Cee[-3]).
\end{equation}

\ssec{The modules of tensor fields}\label{SS:2.7} Let $\fg=\fvect(m|n)$ and $\fg_{\geq}=\mathop{\oplus}_{i\geq
0}\fg_{i}$ in the standard grading ($\deg x_i=1$ for all $i= 1, \dots m+n$). For any other $\Zee$-graded vectorial Lie superalgebra
for whose component $\fg_0$ the notion of lowest weight can be
defined, the construction is identical.

Let $V$ be the $\fg_0=\fgl(m|n)$-module with the \textit{lowest} weight
$\lambda=\lwt(V)$. Let us make $V$ into a~$\fg_{\geq}$-module by setting
$\fg_{+}\cdot V=0$ for $\fg_{+}=\mathop{\oplus}_{i>
0}\fg_{i}$. Let us realize $\fg$ by vector fields on the linear
supermanifold $\cC^{m|n}$ with coordinates $x=(u, \xi)$. The
superspace $T(V)=\Hom_{U(\fg_{\geq})}(U(\fg), V)$ is isomorphic, due
to the Poincar\'e--Birkhoff--Witt theorem, to ${\Cee}[[x]]\otimes
V$. Its elements have a~natural interpretation as formal
\textit{tensor fields of type} $V$. When $\lambda=(a, \dots , a)$ we
will simply write $T(\vec a)$ instead of $T(\lambda)$. We will
usually consider $\fg$-modules coinduced from irreducible
$\fg_0$-modules. For example, $T(\vec 0)$ is the superspace of functions;
$\Vol(m|n)=T(1, \dots , 1\mid  -1, \dots , -1)$ (the bar separates
the first $m$ (``even'') coordinates of the weight with respect to
the matrix units $E_{ii}$ of $\fgl(m|n)$) is the superspace of
\textit{densities} or
\textit{volume forms}.
We denote the generator of $\Vol(m|n)$, as $\Cee[x]$-module,
corresponding to the ordered set of coordinates $x$ by $\vvol_x$.
The space of $\lambda$-densities is $\Vol^{\lambda}(m|n)=T(\lambda,
\dots , \lambda; -\lambda, \dots , -\lambda)$. In particular,
$\Vol^{\lambda}(m|0)=T(\vec \lambda)$ but
$\Vol^{\lambda}(0|n)=T(\overrightarrow{-\lambda})$.

As $\fvect(m|n)$- and $\fsvect(m|n)$-modules,
$\fvect(m|n)=T(\id_{\fgl(m|n)})$.

\ssec{Deformations of the Buttin superalgebra (after
\cite{LSq})}\label{SS:2.8}  As is clear from the definition of the
Buttin bracket, there is a~regrading (namely, $\fb (n; n)$ given by
$\deg\xi_i=0, \deg q_i=1$ for all $i$) under which $\fb(n)$,
initially of depth 2, takes the form
$\fg=\mathop{\oplus}_{i\geq -1}\fg_i$ with
$\fg_0=\fvect(0|n)$ and $\fg_{-1}\cong \Pi(\Cee[\xi])$. Replace now
the $\fvect(0|n)$-module $\fg_{-1}$ of functions (with inverted
parity) by the rank 1 (over the algebra of functions) module of
$\lambda$-densities (also with inverted parity), i.e., we set
$\fg_{-1}\cong \Pi(\Vol^\lambda(0|n))$, where the
$\fvect(0|n)$-action on the generator $\vvol_\xi^\lambda$ is given
by the formula
\begin{equation}
\label{2.7.1} L_D(\vvol_\xi^\lambda) =\lambda \Div D\cdot
\vvol_\xi^\lambda\; \text{ and }\;
p(\vvol_\xi^\lambda)=\od.
\end{equation}

Define $\fb_{\lambda}(n; n)$ to be the Cartan prolong
\begin{equation}
\label{fb} (\fg_{-1}, \fg_0)_*=(\Pi(\Vol^\lambda(0|n)),
\fvect(0|n))_*.
\end{equation}
Clearly, this is a~deform of $\fb(n; n)$. The collection of these
$\fb_{\lambda}(n; n)$ for all $\lambda$'s is called the \textit{main
deformation}, the other deformations, defined in what follows, will
be called \textit{singular}.

The deform $\fb_{\lambda}(n)$ of $\fb(n)$ is a~regrading of
$\fb_{\lambda}(n; n)$ described as follows. For $\lambda
=\frac{2a}{n(a-b)}$, set
\begin{equation}
\label{2.7.2} \fb_{a, b}(n) =\left\{M_f\in \fm (n)\mid a\; \Div
M_f=(-1)^{p(f)}2(a-bn)\frac{\partial f}{\partial \tau}\right\}. 
\end{equation}
For future use, we denote the operator that singles out
$\fb_{\lambda}(n)$ in $\fm (n)$ as follows:
\begin{equation}
\label{2.7.3} \Div_{\lambda}=(bn-aE)\pder{\tau}-a\Delta\;\text{ for
$\lambda =\frac{2a}{n(a-b)}$ and $\Delta=\mathop{\sum}_{i\leq
n}\frac{\partial^2 }{\partial q_i
\partial\xi_i}$}.
\end{equation}
Taking into account the explicit form of the divergence of $M_{f}$
we get
\begin{equation}
\label{2.7.4}
\begin{array}{ll}
\fb_{a, b}(n) &=\{M_f\in \fm (n)\mid (bn-aE)\pderf{f}{\tau} =
a\Delta
f\}=\\
&\{D\in\fvect(n|n+1) \mid L_{D}(\vvol_{q, \xi,
\tau}^a\alpha_{0}^{a-bn})=0\}.\end{array}
\end{equation}
It is subject to a~direct verification that $\fb_{a, b}(n)\simeq
\fb_\lambda(n)$ for $\lambda =\frac{2a}{n(a-b)}$. This isomorphism
shows that $\lambda$ actually runs over $\Cee P^1$, not $\Cee$ as
one might hastily think. The Lie superalgebras $\fb_\lambda(n)$
are simple for $n>1$ and $\lambda\neq 0$, 1, $\infty$ for reasons
clear from eq.~\eqref{2.7.5}. It is also clear that the
$\fb_{\lambda}(n)$ are non-isomorphic for distinct $\lambda$'s, bar
occasional isomorphisms, see \cite{LSh}.

The Lie superalgebra $\fb(n)=\fb_{0}(n)$ is not simple: it has a $\eps$-dimensional  center. Observe that
$\fb_{1}(n)$ and $\fb_{\infty}(n)$
are not simple either. The corresponding exact
sequences are
\begin{equation}
\label{2.7.5}
\begin{array}{c}
0\tto \Cee M_{1} \tto \fb(n)\tto \fle(n)\tto 0,\\
0\tto \fb_{1} ^{(1)}(n)\tto \fb_{1}(n)\tto \Cee\cdot
M_{\xi_1\dots\xi_n} \tto 0, \ \codim  \fb_{1} ^{(1)}(n)= \eps^{n+1};\\
0\tto \fb_{\infty} ^{(1)}(n)\tto \fb_{\infty}(n)\tto \Cee\cdot
M_{\tau\xi_1\dots\xi_n} \tto 0, \ \codim\fb_{\infty} ^{(1)}(n)= \eps^{n}.\\
\end{array}
\end{equation}
Clearly, at the exceptional values of $\lambda$, i.e., 0, 1, and
$\infty$, the deformations of $\fb_{\lambda}(n)$ should be
investigated extra carefully. As we will see immediately, it pays:
at each of the exceptional points we find extra deformations. An
exceptional deformation at $\lambda=-1$ remains inexplicable. Other
exceptional values ($\lambda =\frac12$ and $-\frac32$) come from the
isomorphisms $\fh_{\lambda}(2|2)\cong \fh_{-1-\lambda}(2|2)$ and $\fb_{1/2}(2; 2)\cong \fh_{1/2}(2|2)=\fh(2|2)$.

\ssec{Deforms of $\fb_{\lambda}(n)$ (\cite{LSq})}\label{thdefb}
\textit{For $\fg=\fb_{\lambda}(n)$, set $H=H^2(\fg; \fg)$ for brevity. Then}

1) \textit{$\sdim ~H=(1|0)$ for $\fg=\fb_{\lambda}(n)$ unless
$\lambda=0$, $-1$, $1$, $\infty$ for $n>2$. For $n=2$, in addition
to the above, $\sdim ~H\neq (1|0)$ at $\lambda=\frac12$ and
$\lambda=-\frac32$.}

2) \textit{At the exceptional values of $\lambda$ listed in heading
$1)$ we have}

$\sdim ~H=(2|0)$ \textit{at $\lambda=\pm 1$ and $n$ odd, or
$\lambda=\infty$ and $n$ even, or $n=2$ and $\lambda=\frac12$ or}
$\lambda=-\frac32$.

$\sdim ~H=(1|1)$ \textit{at $\lambda=0$, or $\lambda=\infty$ and $n$ odd,
or $\lambda=\pm 1$ and $n$ even}.

\noindent \textit{The corresponding cocycles $C$ are given by the following
nonzero values in terms of the generating functions $f$ and $g$,
where $d_{\od}(f)$ is the degree of $f$ with respect to odd
indeterminates only; for $k=(k_{1}, \dots , k_{n})$  we set
$q^k=q_{1}^{k_{1}}\dots q_{n}^{k_{n}}$ and $|k|=\sum k_{i}$}:
\[\tiny
\begin{tabular}{|c|c|c|}
 \hline
$\fb_{\lambda}(n)$&$p(C)$&$C(f, g)$\cr \hline \hline
$\fb_{0}(n)$&\textit{odd}&$(-1)^{p(f)}(d_{\od}(f)-1)(d_{\od}(g)-1)fg$\cr
\hline \hline $\fb_{-1}(n)$&$n+1\pmod 2$&$f=q^k, \ g=q^l \longmapsto
(4-|k|-|l|)q^{k+l}\xi_{1}\dots \xi_{n}+$\cr &&$\tau
\Delta(q^{k+l}\xi_{1}\dots \xi_{n})$\cr \hline
$\fb_{1}(n)$&$n+1\pmod 2$&$f=\xi_{1}\dots \xi_{n}, \ g\longmapsto
\begin{cases}(d_{\od}(g)-1)g&\text{if $g\neq af$, $a\in\Cee$}\\
 2(n-1)f &\text{if $g=f$ and $n$ is even}\end{cases}$
 \cr
\hline $\fb_{\infty}(n)$&$n\pmod 2$&$f=\tau\xi_{1}\dots \xi_{n}, \
g\longmapsto
\begin{cases}(d_{\od}(g)-1)g&\text{if $g\neq af$, $a\in\Cee$}\\
2f &\text{if $g=f$ and $n$ is odd}\end{cases}$\cr \hline \hline
$\fb_{\frac12}(2)$&\textit{even}&\textit{see} \eqref{cocycle},
\eqref{bbdef}\cr \hline
\end{tabular}
\]

\textit{On $\fb_{\frac12}(2)\simeq\fb_{-\frac32}(2)\simeq\fh(2|2; 1)$ \emph{(the latter being the non-standard W-grading in which $\deg\xi=0$, other indeterminates being of degree 1)}, the
cocycle $C$ is the one induced on $\fh(2|2)=\fpo(2|2)/\Cee\, K_1$
by the usual deformation (quantization) of $\fpo(2|2)$: we first
quantize $\fpo(2|2)$ and then take the quotient modulo the center}
(generated by constants).

3) \textit{The space $H$ is diagonalizable with respect to the Cartan
subalgebra of $\fder ~\fg$. Let the cocycle $M$ corresponding to the
main deformation be one of the eigenvectors. Let $C$ be another
eigenvector in $H$, it determines a~singular deformation from
heading $2)$. The only cocycles $kM+lC$, where $k, l\in\Cee$, that
can be extended to a~global deformation are those for $kl=0$, i.e.,
either $M$ or $C$.}

\textit{All the singular deformations of the bracket $\{\cdot,
\cdot\}_{old}$ in $\fb_{\lambda}(n)$, except the ones for
$\lambda=\frac12$ or $-\frac32$ and $n=2$, have the following very
simple form even for $p(\hbar)=\ev$}:
\begin{equation}
\label{cocycle} \{f, g\}_{\hbar}^{sing}=\{f, g\}_{old}+\hbar\cdot
C(f, g)\textit{ for any } f, g\in \fb_{\lambda}(n).
\end{equation}

Since the elements of $\fb_{\lambda}(n)$ are encoded by functions
(for us: polynomials or formal power series) in $\tau$, $q$ and
$\xi$ subject to one relation with an odd left hand side in which
$\tau$ enters, it seems plausible that the bracket in
$\fb_{\lambda}(n)$ can be, at least for generic values of parameter
$\lambda$, expressed solely in terms of $q$ and $\xi$. This is,
indeed, the case, and here is an explicit formula (in which $\{f,
g\}_{B.b.}$ is the usual antibracket and
$\Delta=\mathop{\sum}_{i\leq n}\frac{\partial^2 }{\partial
q_i
\partial\xi_i}$):
\begin{equation}
\label{bbdef}
\renewcommand{\arraystretch}{1.4}
\begin{array}{ll}
  \{f_1, f_2\}_{\lambda}^{main}=&\{f_1,
f_2\}_{B.b.}+\lambda(c_{\lambda}(f_1, f_2)f_1\Delta f_2 +
(-1)^{p(f_1)}c_{\lambda}(f_2, f_1)(\Delta f_1)f_2),\\
&\text{where $ \displaystyle c_{\lambda}(f_1, f_2)=\frac{\deg
f_1-2}{2+\lambda(\deg f_2 -n)} $},
\end{array}
\end{equation}
and $\deg$ is computed with respect to the standard grading $\deg
q_{i}=\deg \xi_{i}=1$.


\ssec{Vectorial Lie (super)algebras for $p>0$}\label{p>0} Recall that we are interested in simple Lie (super)algebras; their central extensions and the algebras of derivations, though no less interesting per se, are next objects on our agenda. Every simple $\Zee$-graded Lie (super)algebra $\fg=\oplus \fg_i$ is \textit{transitive}, i.e., such that
\[
[\fg_i, x]=0\text{~~for a given $x\in \fg_j$, where $j>0$, and any $i<0$ implies $x=0$}.
\]
Over any field $\mathbb{K}$ of characteristic $p>0$, in order for the analogs of vectorial Lie (super)algebras be \textit{transitive}, we must change the definition in the two places:

(A) Consider not polynomial coefficients but divided powers in $m$ even and $a-m$ odd indeterminates, whose powers are bounded by the shearing vector $\tilde\un = (N_1,..., N_m, 1, \dots, 1)$, usually abbreviated to $\un = (N_1,..., N_m)$, forming the supercommutative superalgebra (here
$p^{\infty}:=\infty$)
\begin{equation*}\footnotesize
\label{u;N} \cO(a; \tilde\un):= \cO(m; \un|n):=\mathbb{K}[u;
\un]:=\text{Span}_{\mathbb{K}}\left(u^{(\underline{r})}\mid 0\leq r_i
\begin{cases}< p^{N_{i}}&\text{for $i\leq m$}\\
\leq 1&\text{for $m<i\leq a$}\end{cases}\right),
\end{equation*}
 where $u^{(\underline{r})}=\prod u_i^{(\underline{r_i})}$. Set $\One:=(1,\dots, 1)$.

(B) Introduce
\textit{distinguished} partial derivatives $\partial_i$ each of them serving as several
partial derivatives at once, for each of the generators $u_i$, $u_i^{(p)}$, $u_i^{(p^2)}$, \dots (or, in terms of $y_{i,j}:=u_i^{(p^{j-1})}$):
\[\footnotesize{
\partial_i(u_j^{(k)}):=\delta_{ij}u_j^{(k-1)}}\ \text{ for all $k$, i.e., }\partial_i=\sum_{j\geq 1} (-1)^{j-1}y_{i,1}^{p-1}\dots y_{i,j-1}^{p-1}\partial_{y_{i,j}}.
\]
The (general) Lie (super)algebra of vector fields is
\[
\fvect(m;\un|n):=\left\{\sum f_i\partial_i\mid f_i\in \cO(m; \un|n)\right\}.
\]
\textbf{In what follows, speaking about Lie SUPERalgebras, we assume that $p>2$} unless otherwise specified.

\sssec{Linguistics: names. Priorities} In the old literature, the Lie algebra $\fvect(m;\un)$ of vector fields with coefficients in the algebra of divided powers $\cO(m; \un)$  was called \lq\lq the general Lie algebra of Cartan type";  in the modern literature, it is called the \textit{Jacobson--Witt algebra}; it is usually denoted by $W(m;\un)$ for any $m$ and $\un$ whereas the name \textit{Witt algebra} is reserved to the particular case $\fvect(1;\underline{1})$. The more general case $\fvect(1;\un)$ is called \textit{Zassenhaus} algebra if $\un\neq \One$. Jacobson--Witt algebras are simple if $m>1$; there is no special name for the simple derived algebra $\fvect^{(1)}(1;\underline{n})$. We denote the Lie algebra of \textit{divergence free} (or \lq\lq special") vector fields  $\fsvect(m;\un)$, usually abbreviated to $S(m;\un)$.

Because $\fk(1;\un|0)\simeq\fvect(1; \un|0)$, we may assume  that $n>0$ speaking about the contact Lie algebras $\fk(2n+1;\un|0)$. In what follows, for $p=0$, just ignore the shearing vector $\un$. For the proof on non-existence of NIS on $\fvect(m;\un)$, see \cite{Zus}.

Dzhumadildaev generalized a result due to Block (1958) and listed simple $\Zee$-graded vectorial Lie algebras of the four series of \lq\lq Cartan type" in characteristic $p\geq 3$ that have NIS, see \cite{Dz1}. Proofs appeared later in \cite[\S2]{Dz2}. (Dzhumadildaev also did the same, albeit in a somewhat implicit way and also, as  in \cite{Dz1}, without proofs, for infinite-dimensional vectorial Lie algebras over $\Cee$ in \cite{Dz3}.
Later, Farnsteiner obtained the same results in the modular case, see \cite{Fa}.) Dzhumadildaev's and Farnsteiner's proofs differ in details: Dzhumadildaev described explicitly the coadjoint module of the corresponding Lie algebras, while Farnsteiner computed in the universal enveloping algebra. The latter proof is reproduced in the book \cite{SF}.

To extend Dzumadildaev's result to filtered deforms of the algebras he considered, we reproduce an explicit description of these NISes  from \cite{SF}, and of deforms, see \cite{BGLd} and \cite{DK}.

\ssec{Filtered deforms of $\fsvect$ and $\fh$} The Lie (super)algebras of series $\fsvect$ and $\fh$ are not simple, but several of their simple relatives are subquotients of their filtered deforms.

\sssec{Types of  Lie algebras $\fsvect$ described by Tyurin and Wilson \cite{Tyu,W}}
Let $\widehat \cO(m)$ be the algebra of formal power series in $u$ for $\un=(\infty, \dots, \infty)$.  For $p>3$, suppose that
\be\label{distN}
N_1=\dots=N_{t_1}<N_{t_1+1}=\dots=N_{t_2}<\dots<N_{t_{s-1}+1}=\dots=N_{t_s}=N_m.
\ee
Wilson proved that there are only the following three types of non-equivalent classes of volume forms, and hence filtered deforms of divergence-free algebras preserving them:
\begin{equation}\label{svectH}
\renewcommand{\arraystretch}{1.4}
\begin{array}{l}
\fsvect_{h}(m; \un):=\{D\in\fvect(m; \un)\mid
L_D(h\vvol_u)=0\},\text{~~where $h$ is one of the following:}\\
h=\begin{cases}1,\\
1+\bar u, \text{~~where $\bar u:=\prod \bar u_i$ and $\bar
u_i:=u_i^{(p^{\un_i}-1)}$,}\\
\exp(u_{t_i}^{(p^{\un_{t_i} })}):=\mathop{\sum}\limits_{j\geq
0}(u_{t_i}^{(p^{\un_{t_i}}) })^{(j)}.\\
\end{cases}
\end{array}
\end{equation}
For brevity, set $\fsvect_{\exp_i}(m; \un):=\fsvect_{\exp\left(u_{t_i}^{\left(p^{\un_{t_i} }\right)}\right)}(m; \un)$  and $\fsvect:=\fsvect_1$.

\paragraph{Remarks} 1) Kac observed that if $L_D(\exp(u_{t_i}^{(p^{\un_{t_i} })})\vvol_u(\un))=0$, then $D\in\fvect(m;\un)$ with all coordinates of $\un$ finite although $\exp(u_{t_i}^{(p^{\un_{t_i} })})\in\widehat \cO(1)\subset\widehat \cO(m)$, see  \cite{KfiD}.

2) For $p=3$ and 2, the  deforms \eqref{svectH} of $\fsvect$ are also possible; nobody knows if there are no other, non-isomorphic, ones, whereas for $p=2$ there definitely is at least one more deform, its existence is the most spectacular result of \cite{BGLLS}.

3) S.~Tyurin \cite{Tyu} described the Lie algebras of divergence-free type and got an extra type of volume forms, as compared with Wilson's list \eqref{svectH}; Tyurin missed an equivalence.

3) S. Kirillov \cite{Kir} verified Skryabin's remark in \cite{Sk2}, namely for which $i$ is the $i$th derived algebra  from Wilson's list~\eqref{svectH} simple and what  is its dimension:\be\label{dimSvect}
\begin{array}{ll}
\dim\fsvect_{\exp_j}(m; \un)=(m-1)p^{|\un|}&i=0;\\
\dim\fsvect_{1+\bar u}^{(1)}(m; \un)=\dim\fsvect^{(1)}(m; \un)=(m-1)(p^{|\un|}-1)&i=1.
\end{array}
\ee

\sssec{Types of Hamiltonian Lie algebras classified by S.~Skryabin \cite{Sk1,Sk2}} Let $\fh(2k;\un)$ be the $\Zee$-graded Lie algebra  preserving the symplectic form
$\omega_0 = \mathop{\sum}_{1\leq i\leq k}d u_i\wedge du_{k+i}$. Its nonisomorphic filtered deforms are only the following $\fh_{\omega_{i}}(2k;\un)$, where $i=1,2$, preserving the respective forms.
For $k>1$, let $\fh_{\exp_j}(2k;\un)$ be the Lie algebra  that preserves the form (Skryabin calls it \textit{of $2$nd type})
\be\label{2ndtype}
\omega_{2,j} = d\left(\exp (\eps u_j )\mathop{\sum}\limits_{1\leq i\leq k}u_i du_{k+i}\right), \text{~~where $j = t_1, t_2, \dots , t_s$, see~\eqref{distN}, and $\eps\in\Kee^\times$}.
\ee

S.~Skryabin proved (\cite{Sk2}) that for $p>2$ the only other, inequivalent to \eqref{2ndtype} and to each other indecomposable symplectic forms (Skryabin calls them \textit{of $1$st type}) are those that can be reduced to the following normal shapes
\begin{equation}\label{ham}
\begin{array}{l}
\omega_1 = \omega_0 +  \eps \mathop{\sum}\limits_{1\leq i,j\leq 2k} A_{i,j}d(\bar u_i )\wedge d(\bar u_{j}),\text{~~where $\bar
u_i:=u_i^{(p^{N_i})}$, \textbf{NOT} $\bar
u_i:=u_i^{(p^{N_i}-1)}$,}\end{array}
\end{equation}
and $A=(A_{i,j})$ can only be equal to one of the following:
\begin{equation*}\label{3.3}
\renewcommand{\arraystretch}{1.4}
\begin{tabular}{|l|l|l|}\hline
type of $A$&shape of $A$&detailed notation of $\omega_1$\\
\hline
$J_k(0)$&$\antidiag(J_k(0),
-J_k(0)^T)$&$\omega_{1,0}$ for $k>1$\\
\hline
$J_{k,r}(\lambda)$, where $\lambda\neq 0$&$\antidiag(J_{k,r}(\lambda), -(J_{k,r}(\lambda))^T)$&$\omega_{1,r,\lambda}$ for $k=r m$ for $r,m\geq 1$\\
\hline
$C_{k}$&$\antidiag(C_{k}, -C_{k}^T)$&$\omega_{1,C}$ for $k>1$ \\
\hline
\end{tabular}
\end{equation*}
where $J_k(\lambda)$ is a Jordan $k\times k$ block with eigenvalue $\lambda$, and $J_{k,r}(\lambda)$ is a $k\times k$ block matrix with blocks of size $r\times r$, so $k=r\times m$ for some $r,m\geq 1$:\begin{equation*}\label{zhut}
J_{k,r}(\lambda)=\begin{pmatrix}0_r&1_r&&0\\
&&\ddots&\\
0&0& &1_r\\
J_r(\lambda)&0&&0_r\\
\end{pmatrix},\text{~ and $C_{k}=\begin{pmatrix}0&1&&0\\
&&\ddots&\\
0&0& &1\\
1&0&&0\\
\end{pmatrix}$ is of size $k\times k$ for $k>1$.}
\end{equation*}

\sssec{The two conditions on $J_{k,r}(\lambda)$ and $C_k$}\label{2cond} 1) The  case with $J_{k,r}(\lambda)$  occurs only when
\be\label{=}
N_1 + ... + N_{nr} = N_{nr+1} + ... + N_{2nr} \text{~~(recall that $k=r n$)}
\ee
and, furthermore,    $N_{ir - j } = N_{ir}$  for all $i = 1,..., 2n$  and
all  $j = 1,..., r-1$, i.e., the $r$ indeterminates in
each of the $2n$ successive groups have equal heights.

The case with $C_k$ occurs only when condition \eqref{=} is violated.

2) Let $G$ be  the group generated by the cyclic
permutations of the coordinates of vectors in  $\Nee^k$. The 2nd condition requires the identity element to be the only element in  $G$  that fixes the two vectors $v=(N_1,...,N_k)$  and  $w=(N_{k+1},...,N_{2k})$   simultaneously. It suffices to consider representatives of equivalence classes of pairs $(v,w)$ under the $G$-action.

\paragraph{Implicit brackets} For a basis in~$\fh_\omega(2k;\un)$ we take Hamiltonian vector fields
\be\label{skobka}
  X^{\omega}_H = \sum_{i,j} (\partial_i H) (\omega^{-1})^{ij} \partial_j \text{~~or $[X^{\omega}_F,\ X^{\omega}_G]=X^{\omega}_{\{F, G\}_\omega}$},
\ee
where the generating functions run over the union of a basis in the maximal ideal of~$\cO(2k;\un)$ and the collection of non-existing in $\cO(2k;\un)$ for the vector $\un$ with finite
coordinates
Hamilto\-ni\-ans~$u_i^{(p^{\un_i})}$ for $i=1,\ldots,2k$. These  Hamilto\-ni\-ans, though non-existing, generate the vector
fields~$X^{\omega}_{u_i^{(p^{\un_i})}}\in\fvect(2k;\un)$. We identify the basis
element~$X^{\omega}_H$ with its Hamiltonian~$H$.  
The shape of the bracket $\{F, G\}_\omega$,
given in \cite[Lemma 8]{Kir} is still rather implicit, unlike that of $\omega$.
Regrettably, the explicit expressions of $\{F, G\}_\omega$ are given at the moment only in the Ph.D. thesis of S.~Kirillov written in Russian and unpublished.

\ssec{How conditions for simplicity change under the passage from $p=0$ to $p>0$} The \lq\lq natural" objects are vectorial Lie algebras obtained as a result of a (generalized) Cartan prolongation. Such objects are often not simple, the simple \lq\lq derived" (figuratively speaking) of these objects are their first or second derived (commutant) algebras, perhaps, modulo center.

$\bullet$ S.~Kirillov \cite{Kir} checked for which $i$ the $i$th derived algebra of the Hamiltonian Lie algebra from Skryabin's list  is simple and what is its dimension:
\be\label{dimHam}
\begin{array}{l}
\dim\fh_\omega^{(i)}(2k;\un)=\left\{\begin{array}{lll}p^{|\un|}-2&\text{if~}\omega=\omega_0&\text{then~}i=2,\\
 p^{|\un|}&\text{if~}\omega=\omega_2 \text{~~for $k+1\not\equiv 0\pmod p$}&\text{then~}i=0,\\
 p^{|\un|}-1&\text{if~}\omega=\omega_2 \text{~~for $k+1\equiv 0\pmod p$}&\text{then~}i=1,\\
 p^{|\un|}-1&\text{if~}\omega=\omega_1 \text{~~for $\det A\neq0$}&\text{then~}i=1,\\
p^{|\un|}-2&\text{if~}\omega=\omega_1 \text{~~for $\det A=0$ (type $J_s(0)$)}&\text{then~}i=2.\end{array}\right.
\end{array}
\ee

$\bullet$  Over $\Cee$, the supervarieties of parameters of deformations of Poisson and Hamiltonian Lie superalgebras can differ, see \cite{LSq}. For $p=2$, the forms $\omega_1$, see eq. \eqref{ham}, do not exist; but instead there is a 1-parametric family of non-isomorphic deforms different from the above --- desuperisations of the Lie superalgbra $\fb_{a:b}(n;\un)$, see \cite{BGLLS}.

\begin{equation}\label{simplCond}
\begin{tabular}{|l| l|}
\hline
$\fg$&conditions for simplicity\\
\hline
\hline
$\fvect(m;\un|n)$&unless $m=1, \ n=0$ or $m=0,\ n<2$\\
$\fvect^{(1)}(1;\un|n)$&$\un>2$\\
\hline
\hline
$\fsvect(m;\un|n)$, $\fsvect_{(1+\bar u)}(0|n)$&if $m=0$, then $n>2$\\
$\fsvect^{(1)}(m;\un|n)$&$(m, n)\neq(1,0)$\\
\hline
$\fsvect_{\exp_i}(m;\un|n)$&$(m, n)\neq(1,0)$\\
\hline
$\fsvect_{(1+\bar u)}^{(1)}(m;\un|n)$&$(m, n)\neq(1,0)$\\
\hline
\hline
$\fh_{\omega_0}^{(2)}(2k;\un)$&\\
$\fh_{\omega_{1,j}}(2k;\un)$,  where $k+1\not\equiv 0\pmod p$&\\
$\fh_{\omega_{1,j}}^{(1)}(2k;\un)$,  where $k+1\equiv 0\pmod p$&\\
$\fh_{\omega_2}^{(1)}(2k;\un)$,  where $\det A\neq0$&\\
$\fh_{\omega_2}^{(2)}(2k;\un)$,  where $\det A=0$ (type $J_s(0)$)&\\
\hline\hline

$\fk(2k+1;\un)$ if $2k+4\not\equiv 0\pmod p$&always, $p>2$\\
$\fk^{(1)}(2k+1;\un)$ if $2k+4\equiv 0\pmod p$&always, $p>2$\\
\hline

\end{tabular}
\ee

\section{On existence of NIS on vectorial Lie (super)algebras}\label{sources} We see two types of reasons why there is a NIS on a given $\fg$; let us consider them in detail.

(1) \textbf{There is only one geometric structure leading to NIS: volume}. Invariant with respect to any changes of indeterminates (coordinates) unary operators between sections of tensor fields on a (super)manifold $M$ with fibers irreducible under the action of the Lie (super)al\-gebra of linear changes of indeterminates are classified for the case where functions are polynomi\-als or power series over $\Cee$, and, for $M=S^1$ (i.e., when functions considered are Laurent polynomials or series) and over fields of characteristic $p$, see \cite{BL} and references therein. The answer is as follows: if the fibers are finite-dimensional, the only operators  are the exterior differential,  the (super)residue (when  $M$ is the complexified supercircle  associated with the trivial bundle or the Whitney sum of the trivial and M\"obius bundle) or the Berezin integral or its version for $p>0$, i.e., the coefficient of the highest possible degree in the Taylor series decomposition of $f(u)$ for $f(u)\vvol_u$.

(2) \textbf{An algebraic reason for existence of a NIS} on the $\Zee$-graded Lie (super)algebras $\fg$ is a necessary condition, see Lemma~\ref{Pr 4.7}.

\ssec{Lemma (Corollary 4.9 in Ch.3 in \cite{SF})}\label{Pr 4.7} \textit{Let $\fg=\mathop{\oplus}_{-d\le i \leq h}\fg_i$ be a finite-dimensional simple $\Zee$-graded Lie algebra, $B$ a NIS on $\fg$. Then}
\[
\begin{array}{ll}
B(\fg_i, \fg_j)=0&\text{ for $i+j\neq h-d$};\\
\dim\fg_i=\dim\fg_{h-d-i}.
\end{array}
\]

From Lemma~\ref{Pr 4.7} one deduces --- with certain effort ---  the following statement.

\ssec{Corollary (Theorems 6.3 and 6.4 in Ch. 4 in \cite{SF})}\label{Th6.4}  1)
\textit{The Lie algebra $\fvect(n;\un)$ has a NIS if and only if either $n=1$ and $p=3$ in which case  NIS is
\be\label{bilFvect1}
(u^{(a)}\partial, \ u^{(b)}\partial):=\int u^{(a)}u^{(b)} du\pmod 3;
\ee
or $n=p=2$ in which case NIS is
\be\label{bilFvect2}
(u^{(a)}\partial_i, \ u^{(b)}\partial_j):=(i+j)\int u^{(a)}u^{(b)} du_1\wedge du_{2}\pmod 2,
\ee
where $\int f(u) du_1\wedge\cdots\wedge du_{n}:=\text{coefficient of~}u^{(\tau(\un))}$ and $\tau(\un):=(p^{N_1}-1, \dots, p^{N_n}-1)$.}

2) \textit{If $n>2$, then $\fsvect^{(1)}(n;\un)$ has a NIS if and only if $n=3$; explicitly
\begin{equation}\label{nisSvect}
\begin{array}{l}
(\partial_i, \ D_{jk}(u^{(\tau(\un))}))=\sign(i,j,k),\\
\end{array}
\end{equation}
where \emph{(for possible future use we give the formula for the general super case)}
\begin{equation*}\label{basisSvect}
D_{ij}(f):=\partial_j(f)\partial_i - (-1)^{p(f)(p(\partial_i)+p(\partial_j))+p(\partial_i)p(\partial_j)}\partial_i(f)\partial_j\text{~~for any $f\in\cO(m;\un|n)$},
\end{equation*}
and extending the form \eqref{nisSvect} to other pairs of elements by invariance.}


\ssec{Series $\fvect$} Classification allows us to forget about deforms of $\fvect(m;\un)$ for $p>3$ and study other series. For example,  the  true deforms of Lie algebras $\fvect(1;\un)$ are described in \cite{DK}; they are either certain filtered deforms, or isomorphic to certain Hamiltonian Lie algebras, see the next Subsection, where $\fh_{\omega}^{(\infty)}$ designates the appropriate simple derived algebra (1st or 2nd).

\ssec{Series $\fsvect$}  Due to the following \lq\lq occasional isomorphisms" (taking place for reasons given in parentheses below) we do not consider $\fsvect$ for the following values of parameters:
\begin{equation*}\label{IsoSvect}
\begin{array}{ll}
\fsvect(2;\un|0)\simeq \fh(2;\un|0)\text{~~ for any $p$}&\text{(obviously)},\\
\fsvect^{(1)}(3;\One|0)\simeq \fh^{(1)}(4;\One|0)\simeq\fp\fs\fl(4)\text{~~ for $p=2$}&\text{(as shown in \cite{ChKu})}.\\
\end{array}
\end{equation*}

\sssec{Conjecture}[Proved for $n=p=3$ and $\un=(111)$]\label{ClSvect} \textit{There is no NIS on simple Lie algebras} $\fsvect_{\exp_i}(n; \un)$ and $\fsvect_{1+\bar u}^{(1)}(n; \un)$.

\ssec{Series $\fh$} For $p>2$, consider the Lie algebra $\fh(2k; \un)$ and a standard symplectic form $\omega_0 =
\sum_{1\leq i\leq k} du_i\wedge du_{k+i}$, set
\[
  (F,\ G)_{\omega_0} := \int F G\ (\omega_0)^{\wedge k} = \int F G \ \,du_1\wedge\cdots\wedge du_{2k}\text{~~for any~~}F,G\in\fh(2k; \un).
\]
\sssec{Lemma}\label{hNIS} \textit{For all symplectic forms~$\omega$ given by the expressions~\eqref{2ndtype} and \eqref{ham}, a NIS
on~$\fh^{(\infty)}_\omega(2k;\un)$ is given by the following formula}:
\be\label{NISHam}
  (X^\omega_F,\ X^\omega_G) := (F,\ G)_\omega = \int F G\,(\omega)^{\wedge k} .
\ee

\begin{proof}
Since $\{H,\cdot\}_{\omega} = X^{\omega}_H$ is a derivation, the invariance of~\eqref{NISHam} can be written as
 \[
    (\{F, G\}_\omega,\ H)_\omega - (F,\ \{G,H\}_\omega)_\omega
    = \int \left(\{F,G\}_\omega H - F\{G,H\}_\omega \right)\,(\omega)^{\wedge k}
    = \int \{G, FH\}_\omega\,(\omega)^{\wedge k}=0.
 \]
It is easy to check that
\[
 \int \{F, G\}_\omega\,(\omega)^{\wedge k}=0\text{~~ for all
$F,G\in\fh^{\infty}(2k;\un)$.}
\]
 Observe that if $X^{\omega_0}_F, X^{\omega_0}_G\in\fh_{\omega_0}^{(\infty)}(2k;\un)$, then $X^{\omega}_F,X^{\omega}_G\in\fh_{\omega}^{(\infty)}(2k;\un)$ for all symplectic forms given by formulas~\eqref{2ndtype} and \eqref{ham}.
We see that $(\omega)^{\wedge k} = (\omega_0)^{\wedge k} + O(\eps)$, where $O(\eps)$ is the term proportional to $\eps$.
If $(F,\ G)_{\omega_0} \neq 0$ for some non-zero $X^{\omega_0}_F,X^{\omega_0}_G\in\fh_{\omega_0}^{(\infty)}(2k;\un)$, then
 \[
 (F,\ G)_\omega = (F,\ G)_{\omega_0} + O(\eps).
 \]
Since $(F,\ G)_{\omega_0}\neq 0$ and $\eps\in\Kee$ is arbitrary, we see that
$(F,\ G)_\omega \neq 0$. Therefore, the form~$(\cdot,\ \cdot)_\omega$ is non-degenerate on~$\fh_\omega^{(\infty)}(2k;\un)$.
\end{proof}

\ssec{Series $\fk$ and its subalgebras}\label{sssK}  On the supermanifold $\cK^{2n+1|m}$ with a contact structure given by $\alpha:=dt+\sum p_idq_i+\sum \theta_jd\theta_j$, consider the space of $\lambda$-densities, 
\[
\cF_{\lambda}(\un):=
\begin{cases}\cF\alpha^{\lambda/2}&\text{if $(m,n)\neq (0,0)$},\\
\cF\alpha^{\lambda}&\text{if $(m,n)=(0,0)$}, 
\end{cases}
\]
where $\cF$ is the space of \lq\lq functions" and $\lambda\in\Kee$. We have $\Vol(\un)\simeq\cF_{2n+2-m}(\un)$, where  $\Vol(\un)$ is the space of volume forms. For $p\neq 2$,  NIS obviously exists on the space of $\nfrac12$-densities, i.e., on $\cF_{n+1-k}(\un)$,  where $m=2k$. 

Observe that this NIS can be odd, but there are two natural ways to define the parity of it, cf. \cite{MaG}; in our situation, where the integral is considered to be even for any $n$ if $m=0$, the adequate  definition is $p(\text{NIS})=m\pmod2$:
\[
p(\text{NIS})= \begin{cases}2n+1-m\pmod 2&\\
m\pmod2.&\end{cases}
\]
Since 
\[
\fk(2n+1;\un|m)\stackrel{\text{as $\fk(2n+1;\un|m)$-modules}}{\simeq}
\begin{cases} \cF_{-2}(\un)&\text{if $(m,n)\neq (0,0)$},\\
\cF\cF_{-1}(\un)&\text{if $(m,n)=(0,0)$}, 
\end{cases}
\]  
it follows that
\be\label{cond}
\fk(2n+1;\un|m)\text{\ has NIS $(f\sqrt{\vvol}, g\sqrt{\vvol}):=\int fg\vvol$ if~~}2n+2-m\equiv -4\pmod p.
\ee

\ssec{No NIS on $\fk\fa\fs$ and $\fk\fle$} Eq. \eqref{cond} shows there is a NIS on $\fk(1;\un|6)$, so we have to verify if its restriction to $\fk\fa\fs$, see eq.~\eqref{nonpos}, is non-degenerate. Since $K_1\in\fk\fa\fs$, the only element it can be paired with a non-zero result is the highest possible power of all indeterminates but $t^k\xi_1\xi_2\xi_3\eta_1\eta_2\eta_3\not\in\fkas$ for any $k$ and any $p$, so no NIS.

The Lie superalgebra $\fk\fle(9;\un|6)$ has no  NIS by argument \eqref{cond} and thanks to Lemma~\ref{Pr 4.7}.

\ssec{$p=5$: NIS on the Melikyan algebra} This fact is known, see, e.g., \cite{S}. The condition in \eqref{cond} is satisfied for $n=2$, $m=0$. For descriptions of the Melikyan algebra $\fm\fe(\un)$ as a subalgebra in $\fk(5;\un)$,  see \cite{GL}. Therefore, to prove the existence of a NIS on $\fk(5;\un)$,  it suffices, thanks to simplicity of $\fm\fe(\un)$, to verify that the restriction of NIS from $\fk(5;\un)$ onto  $\fm\fe(\un)$ does not vanish on just one pair of elements. (Thanks to Dzhumadildaev's interpretation of $\fm\fe(\un)$ as a deform of the Poisson Lie algebra, see \cite{KD}, this result contributes to Fact~\ref{Fact}.)

\ssec{$p=3$: NIS on some of Skryabin algebras. No NIS on Ermolaev and Frank algebras}\label{Skr} The exceptional simple $\Zee$-graded vectorial Lie algebras known today to be indigenous to  $p=3$ are described in \cite{GL}, where the number $\Par\un$ of parameters, the shearing vector depends on, is correctly computed for the first time for some of these algebras:
\begin{equation}\label{standAlgs}
\begin{array}{l}
\fby^{(1)}(7;\un),\text{~~(the derived of) the ``big" Skryabin algebra, }\\
\text{its divergence-free subalgebra~~}\fsby^{(1)}(7;\un),\\
\fdy^{(1)}(10;\un), \ \  \text{(the derived of) the ``deep" Skryabin algebra, initially denoted by $Z$},\\
\fmy(6; \un), \ \fsmy^{(1)}(6; \un), \ \  \text{the ``middle"Skryabin algebra, initially denoted by $Y$},\\
\fer^{(1)}(3;\un), \ \ \text{(the derived of) the Ermolaev algebra a.k.a. $R$},\\
\ffr^{(1)}(3;\un), \ \ \text{(the derived of) the Frank algebra}.
\end{array}
\end{equation}
Lemma~\ref{Pr 4.7} implies that among the algebras \eqref{standAlgs}, only the following simple Skryabin Lie algebras may have NIS, and they do have it (we proved this for $\un=\One$ only):
\be\label{SkNis}
\fby^{(1)}(7;\un) ,\ \ \fmy(6; \un), \ \ \fs\fmy^{(1)}(6; \un).
\ee
Observe that all of the Skryabin algebras are results of generalized, perhaps partial, Cartan prolongation of the non-positive part of $\fbr(3)$ in one of the $\Zee$-gradings of $\fbr(3)$; both of its Cartan matrices were discovered by Skryabin (for a method proving that these are the only possible Cartan matrices  of $\fbr(3)$, see \cite{BGL}).

\ssec{NIS on the two simple stringy Lie (super)algebras over $\Cee$} For the classification of simple stringy Lie superalgebras over $\Cee$, i.e., vectorial ones on a supercircle, see \cite{GLS}. We see that $\Vol^L(1|n)\simeq\Pi^n(\cF_{2-n}^L)$, where superscript $L$ is for Laurent, and hence the Neveu--Schwarz type stringy Lie superalgebra $\fk^L(1|n)$ has NIS if $2-n=-4$, i.e., $n=6$; this NIS is even.

For the Ramond-type stringy Lie superalgebra $\fk^M(1|n)$ preserving the distribution given by the form equivalent to $dt+ td\theta +\mathop{\sum}_{1\leq i\leq n-1}\theta_id\theta_i$ (so the weight of $\theta$ is equal to 0) on the supercircle associated with the Whitney sum of the $(n-1)$-dimensional cylinder and the Moebius bundle, for $n>1$, we have
\[
\Vol^M(1|n)\simeq\Pi^n(\cF_{2-(n-1)}^L).
\]
Therefore,  $\fk^M(1|n)$ has NIS if $3-n=-4$, i.e., $n=7$; this NIS is an odd one. (For $n=1$, the condition \eqref{cond} turns into $1=-2$. Hence, no NIS.)

In terms of tensor fields, we see that
\[
\text{$\fvect^L(1|N):=T^L(\id_{\fgl(1|N)})$, so $(\fvect^L(1|N))^*:= (T^L(\id_{\fgl(1|N)}))^*= T^L((\id_{\fgl(1|N)})^*\otimes \str)$},
\]
where $\str$ is a 1-dimensional module given by the supertrace, hence $\fvect^L(1|N)\not\simeq (\fvect^L(1|N))^*$ always.

The cases $\fg=\fsvect^L(1|N)$, and its deform $\fsvect_\lambda^L(1|N)$ preserving $t^\lambda\vvol_{(t|\theta)}$, are similar to the cases of $\fvect$, except for $N=2$ with its  decomposition $\fsvect_\lambda^L(1|2)=\oplus_{-1\leq i\leq 1}\fg_i$ where $\fg_0=\fgl(1|1)^{\ell(1)}$ and $\fg_{\pm 1}\simeq V^{\ell(1)}t^{\pm \lambda}$ for the tautological $\fgl(1|1)$-module $V$ and the $V$-valued space of loops $V^{\ell(1)}:= V\otimes \Cee[t^{-1}, t]$;
in this case, a NIS could have existed but does not, see Theorem~\ref{noNIS}.

\section{Odd NIS. Queer Lie superalgebras, queerified Lie algebras, and exceptions}\label{SQueer}

\ssec{Odd NIS on $\fq(n)$}
The \textit{queertrace} $\qtr\colon (X, Y)\longmapsto \tr Y$ on $\fq(n)$, see \eqref{q}, defines a natural NIS on~$\fq(n)$:
\[
  B((X,Y),\ (X^\prime,Y^\prime)) = \qtr( (X,Y)\cdot(X^\prime,Y^\prime)) =  \tr (X Y^\prime + X^\prime Y).
\]

In the same way as the (super)trace  on the associative superalgebra of supermatrices $\Mat(m|n)$ gives rise to a NIS on $\fpsl$, the \textit{queertrace} gives rise to  an odd NIS on the Lie superalgebra $\mathfrak{psq}(n)$ which is simple if $p=0$, and also if $p> 2$ and $n\not\equiv0\pmod p$

\ssec{Queerified Lie algebras}\label{sssQg} The Lie superalgebra $\fq(n)$ is a \lq\lq queerification\rq\rq\ of the associative algebra $\Mat(n)$. If $p=2$, it is possible to \lq\lq queerify\rq\rq\ any Lie algebra, see \cite{BLLSq}.

Let $p=2$ and $\fg$ a \textit{restricted} Lie algebra with a NIS $(\cdot,\cdot)$. Then $\fq(\fg)$ has an odd NIS $q(\cdot,\cdot)$
defined as follows, where $\Pi$ is the change of parity operator:
\be\label{qNIS}
q(x,\Pi(y))= (x,y)\text{~~and $q(x,y)=q(\Pi(x),\Pi(y))=0$ for any $x,y\in \fg$.}\ee

\ssec{Poisson Lie superalgebras} Let us show that both Poisson Lie superalgebras and their particular\footnote{If $p\neq 2$, there are no deforms apart from the quantization; for $p=2$, this is not so, see \cite{BLLS}.} deforms resulting from quantization have NIS if $\un_p=\un_q$, see~\eqref{PoiDefrNpNq}. We do not know if this is so if $\un_p\neq\un_q$.

\sssec{Poisson Lie superalgebras over $\Cee$} On $\fpo(0|m)$ realized on the space $\Cee[\xi, \eta]$, where $\xi=(\xi_1, \dots, \xi_k)$ and  $\eta=(\eta_1, \dots, \eta_k)$ for $m=2k$, or on the space $\Cee[\xi, \eta, \theta]$ for $m=2k+1$, NIS of parity  $p(\int)=m\pmod2$ is defined by the formula
\be\label{PoiNIS}
(f,g)=\int fg\vvol,\text{~~where $\vvol:=\vvol(\xi,\eta)$ or $\vvol(\xi,\eta,\theta)$, }
\ee
and where $\int F\vvol$ is the Berezin integral $=$ the coefficient of the $m$th degree monomial of $F$ in the Taylor series expansion.

Tyutin classified the deforms of $\fpo(0|m)$, see \cite{Ty}. There is just one class $\cQ$ of deformations, called \textit{quantization}; $\cQ$ sends the integral into the super trace (resp. queer trace):
\be\label{PoiDefr}
\cQ: \fpo(0|m)\tto\begin{cases} \fgl(\Cee[\xi])&\text{if $m=2k$}\\
\fq(\Cee[\xi])&\text{if $m=2k-1$}. \end{cases}
\ee

\sssec{NIS on Poisson Lie superalgebras over $\Kee$ for $p>0$}\label{PoQuantiz}
NIS is defined by the direct analog of formula \eqref{PoiNIS}.

\textbf{Aside}. For the case where the shearing vector is of the form $\un=(\un_p,\un_q)$ with $\un_p=\un_q$, and odd indeterminates $\xi=(\xi_1,\dots, \xi_k)$ and $\eta=(\eta_1,\dots, \eta_k)$, the description of quantizations is of the same form as over $\Cee$, i.e.,
\be\label{PoiDefrNpNq}
\cQ: \fpo(2n;\un_p\un_q|m)\tto\begin{cases} \fgl(\cO[q;\un_q\mid \xi])&\text{if $m=2k$}\\
\fq(\cO[q;\un_q\mid \xi])&\text{if $m=2k-1$}. \end{cases}
\ee

There are, however, other deforms of~$\fpo$: e.g., for~$p=5$, there are Melikyan algebras, see~\cite{KD}; for~$p=2$, there is a 1-parametric family --- desuperization of~$\fb_\lambda$, see~\cite{BGLLS}.

\ssec{NIS on \protect{$\mathfrak{sb}(n;\protect\un)$} and \protect{$\mathfrak{sle}^{(1)}(n;\protect\un)$}}\label{ssSB} Recall that by $\fb(n)$ we denote the Lie superalgebra
on the space of functions (divided powers if $p>0$) in $n$ even and $n$ odd indeterminates with Schouten bracket a.k.a. antibracket. Let $\cO(n;\un|n)=\Kee[q;\un, \xi]$ and
\[
\fs\fb(n;\un)\stackrel{\text{as space}}{=}\left\{f\in\cO(n;\un|n)\mid \mathop{\sum}\limits_{i\leq
n}\frac{\partial^2}{\partial_{q_i} \partial_{\xi_i}} (f)=0\right\}.
\]
The direct analog of \eqref{PoiNIS} defines NIS on the Lie superalgebra $\fs\fb(n;\un)$, and its restriction to its subquotient $\fs\fle^{(1)}(n;\un)$ is a NIS. Observe that $p(\text{NIS})\equiv n+1\pmod 2$.

\ssec{Loop superalgebras} For any (simple, or a relative of a simple) finite-dimensional Lie (super)algebra $\fg$ over $\Cee$, we consider spaces $\fg^{\ell(1)}:=\fg\otimes\Cee[t^{-1}, t]$ of \lq\lq loops", i.e., $\fg$-valued functions on the circle  $S$ expandable into Laurent polynomials (or their completions, Laurent series); here $t=\exp(i\varphi)$, where $\varphi$ is the angle parameter on $S$. We introduced notation $\fg^{\ell(1)}$ to distinguish from  $\fg^{(1)}:=[\fg, \fg]$; we give a similar notation for \lq\lq twisted" loop algebras.

Recall that if $\psi$ is an order $m$ automorphism of $\fg$, and $\fg_{\bar k}$, where $0\leq \bar k\leq m-1$, are eigenspaces of $\psi$ with eigenvalue $\exp(\nfrac{2\pi k i}{m})$, where $i=\sqrt{-1}$, then 
\[
\fg_\psi^{\ell(m)}:=\oplus_{0\leq \bar k\leq m-1,\ j\in\Zee}\ \fg_{\bar k}t^{k+mj}. 
\]

Non-isomorphic loop superalgebras $\fg^{\ell(1)}:=\fg\otimes\Cee[t^{-1}, t]$ with values in simple finite-dimensional Lie superalgebras $\fg$ over $\Cee$, and twisted loops $\fg_\psi^{\ell(m)}$ corresponding to order $m$ automorphisms $\psi$ of $\fg$ are classified in~\cite{Se}.

\sssec{NIS on loop superalgebras} The NIS \lq\lq$\tr(x,y)$" for any $x,y\in\fg$ on the target space $\fg$ from the left column in~\eqref{**} is $\tr(\ad_x\ad_y)$ in the purely even case; $\str(\rho(x)\rho(y))$ or $\qtr(\rho(x)\rho(y))$ for an appropriate representation $\rho$ for all superalgebras with an even NIS, except for $\fh^{(1)}(0|m)$ in which case \lq\lq$\tr(H_f,H_g)$" is $\int fg\vvol$. 

\textit{Each of the twisted loop superalgebras $\fg^{\ell(m)}_\psi$ listed in~\eqref{**}, and simple loop algebras listed in~\cite[pp. 54, 55]{K}, has a NIS induced by \lq\lq$\tr$" on $\fg$}:
\be\label{oddNIS}
(a, b)=\Res \lq\lq\tr" (ab)\textit{~~for any $a,b\in \fg^{\ell(m)}_\psi$},
\ee
\textit{where $\Res f(t)$ is the coefficient of $t^{-1}$ in the Laurent polynomial} $f$.

\sssec{Central extensions of loop superalgebras}\label{CentExtLoo} The central extensions of the loop superalgebra with a NIS on the target space $\fg$ is given by the formula
\be\label{LooCenExt}
c(f,g)=\Res \lq\lq\tr" (fdg)\text{~~for any $f,g\in \fg^{\ell(m)}_\psi$},
\ee
where $\Res f(t)dt$) is the coefficient of $t^{-1}dt$ in the  1-form $fdt$. 
Observe that 
\be\label{OddC}
\text{
If \lq\lq$\tr$" is odd, as on $\fpsq(n)$ and $\fh^{(1)}(0|2n+1)$, the center thus obtained is also odd.}
\ee

\paragraph{Loop superalgebras with values in $\fvect(0|n)$ for $n\geq 2$, $\fsvect(0|n)$ and $\widetilde{\fsvect}(0|n)$  for $n\geq 3$, and $\fspe(m)$ for $m\neq 4$}\label{conj1} None of these target algebras has a NIS; the only non-trivial automorphism exists only for $\fvect(0|n)$, see \eqref{*}, but twisted loop superalgebra corresponding to this automorphism is isomorphic to the non-twisted one (\cite{Se}).  Hence, none of  loop superalgebras $\fg^{\ell(1)}$ with values these Lie superalgebras $\fg$  has central extension obtained by recipe~ \eqref{LooCenExt}. \textbf{Conjecturally, these $\fg^{\ell(1)}$ have no central extension at all}.

\paragraph{Loop superalgebras with values in $\fpsl(n|n)$  for $n\geq 2$, $\fpsq(n)$  for $n\geq 3$, $\fh^{(1)}(0|n)$  for $n\geq 4$, and $\fspe(4)$} Each of these target algebras $\fg$ has one non-trivial central extension, whereas $\fpsl(2|2)\simeq\fh^{(1)}(0|4)$ has three of them. Therefore, $\fg^{\ell(1)}$ has infinitely many central extensions, given by loops with values in each of the centers of $\fg$ for each of these target spaces $\fg$.

\sssec{Derivations of loop superalgebras} For $\fg$ simple, the outer derivations the loop superalgebra $\fg^{\ell(1)}$ constitute the semidirect sum $\fG\ltimes M$, where $\fG:=\fder\,\Cee[t^{-1}, t]$ and $M$ is the space of loops with values in $\fout\, \fg:=(\fder\, \fg)/\fg$. The popular affine Kac--Moody algebras correspond (in the non-super case) to just one derivation, $t\nfrac{d}{dt}$.

\ssec{Lie superalgebras of polynomial growth with nonsymmetrizable Cartan matrices}\label{SSnonsymm} Hoyt and Serganova~\cite{HS} proved Serganova's conjecture~(\cite{LSS}) that there are only two types of indecomposable nonsymmetrizable Cartan matrices corresponding to Lie superalgebras of polynomial growth. These Lie superalgebras are the following ones (the algebras of one of these types are double extensions of simple Lie superalgebras, the algebras of the other type are not):

1) The series of twisted loops with values in $\fpsq(n)$ corresponding to the 2nd order automorphism $\Pty(x)=(-1)^{p(x)}x$ for any $x\in\fpsq(n)$. Let $\Pi$ be the change of parity operator; set
\[
\fpsq(n)^{\ell(2)}:=\left(\mathop{\oplus}\limits_{k\in\Zee} \fsl(n)t^{2k}\right)\bigoplus \left(\mathop{\oplus}\limits_{s\in\Zee} \Pi(\fsl(n))t^{2s+1}\right).
\]
A double extension $\widetilde{\fpsq}(n)^{\ell(2)}$ of $\fpsq(n)^{\ell(2)}$ determined by an even central extension and an odd outer derivation  has Cartan matrices with the same Dynkin graphs as that of the loop algebra $\fsl(n)^{\ell(1)}$, but --- in the simplest case --- one of the nodes of the Dynkin graph for $\widetilde{\fpsq}(n)^{\ell(2)}$ being \lq\lq gray", the other Cartan matrices being obtained from the simplest one by \lq\lq odd reflections", see~\cite{CCLL}.

2) The exceptional Lie superalgebra $\widehat\fsvect^L_\alpha(1|2)$ with Cartan matrices \eqref{eq-3d-137} is NOT the double extension of the stringy Lie superalgebra (which is simple for $\alpha\not\in\Zee$)
\begin{align*}
\fsvect_\alpha^L(1|2):={}&\{D\in\fvect^L(1|2)\mid L_D(t^\alpha\vvol(t,\theta_1,\theta_2))=0\}\\
{}={}&\{D=f\partial_t+\sum f_i\partial_{\theta_i}\in\fvect^L(1|2)\mid \alpha f= -t\Div D\}\\
{}={}&(\widehat\fsvect^L_\alpha(1|2))/\fc
\end{align*}
in the sense of definition in Subsection~\ref{DD}: $\widehat\fsvect^L_\alpha(1|2)$ has properties 1) and 2) but not 3); there is no NIS on either $\widehat\fsvect^L_\alpha(1|2)$ or $\fsvect^L_\alpha(1|2)$.

\footnotesize
\sssec{Root system of $\widehat\fsvect^L_\alpha(1|2)$}
For Chevalley generators of the Lie superalgebra we denote $\widehat\fsvect^L_\alpha(1|2)$ take
\begin{equation}
 \label{eq-a2-G1}
\arraycolsep=0pt
\renewcommand{\arraystretch}{1.4}
 \begin{array}{lll}
 X_1^+= \theta_1\delta_{\theta_2},&
 X_2^+=t\delta_{\theta_1},&
 X_3^+=\theta_2t\partial_{t}-(\alpha+1)\theta_1\theta_2\delta_{\theta_1},\\
 X_1^-=\theta_2\delta_{\theta_1},&
 X_2^-=\alpha\fnfrac{\theta_1\theta_2}{t}\delta_{\theta_2}+\theta_1\partial_{t},&
 X_3^-=\delta_{\theta_2},\\
 H_1=\theta_1\delta_{\theta_1}-\theta_2\delta_{\theta_2},\quad &
 H_2=t\partial_{t}+\theta_1\delta_{\theta_1}+\alpha\theta_2\delta_{\theta_2},\quad &
 H_3= t\partial_{t}+(\alpha+1)\theta_1\delta_{\theta_1}.
\end{array}
\end{equation}

By means of the isotropic reflections, see \cite{CCLL}, we establish that  $\widehat\fsvect^L_\alpha(1|2)$ has the following
Cartan matrices:
\begin{equation}
 \label{eq-3d-137}
\begin{pmatrix} 2&-1&-1\\ 1-\alpha&0&\alpha\\ 1+\alpha&-\alpha&0 \end{pmatrix},\quad
\begin{pmatrix} 0&-\alpha+1&-2+\alpha\\ 1-\alpha&0&\alpha\\ -1&-1&2
\end{pmatrix},\quad
\begin{pmatrix} 0&-\alpha&\alpha+1\\ -1&2&-1\\ 1+\alpha&-\alpha-2&0
\end{pmatrix}.
\end{equation}
In order to compare these matrices, reduce them to the following  normal forms having renumbered rows/columns and rescaled
(by definition $\alpha \not\in\Zee$, in particular, $\alpha \neq 0, \pm 1$, so the fractions are well defined), respectively:
\begin{equation}
 \label{eq-3d-138}
\begin{pmatrix} 2&-1&-1\\-1+\frac{1}{\alpha}&0&1\\1+\frac{1}{\alpha}&-1&0
\end{pmatrix},\quad
\begin{pmatrix} 2&-1&-1\\-1+\frac{1}{1-\alpha}&0&1\\1+\frac{1}{1-\alpha}&-1&0
\end{pmatrix},\quad
\begin{pmatrix} 2&-1&-1\\-1+\frac{1}{1+\alpha}&0&1\\1+\frac{1}{1+\alpha}&-1&0
\end{pmatrix}.
\end{equation}
Evidently, the maps $\alpha \mapsto 1+\alpha$ and~$\alpha
\mapsto 1-\alpha$ establish isomorphisms of the Lie superalgebras corresponding to these Cartan matrices, so one may assume that $\RE \alpha
\in\Big(0,\frac{1}{2}\Big]$.

\normalsize
In manual computations, it is convenient  to work  not in the Chevalley basis of $\fsvect_\alpha^L(1|2)$, but in the basis consisting of the following $E_\alpha$-homogeneous elements
\[
\begin{array}{ll}\begin{array}{ll}
L_n&=t^{n}\left(t\partial_t +\nfrac12(n+1+\alpha)(\theta_1\partial_{\theta_1}+\theta_2\partial_{\theta_2})\right),\\

E_n&=t^{n}\theta_2\partial_{\theta_1},\\

F_n&=t^{n}\theta_1\partial_{\theta_2},\\

G_n &=t^{n}\left(\theta_2\partial_{\theta_2}-\theta_1\partial_{\theta_1}\right),\\
\end{array}&
\begin{array}{ll}
\lambda_n&=t^{n-1}\theta_2\left(t\partial_t +(n+\alpha)\theta_1\partial_{\theta_1}\right),\\

\eps_n&=t^{n-1}\theta_1\left(t\partial_t +
(n+\alpha)\theta_2\partial_{\theta_2}\right),\\

\varphi_n&=t^{n+1}\partial_{\theta_1},\\

\gamma_n&=t^{n+1}\partial_{\theta_2}.\\
\end{array}\end{array}
\]
There is only one outer derivation of $\fsvect_\alpha^L(1|2)$ for $\alpha\not\in\Zee$ (for a proof, see \cite{Po}) given (modulo the space of inner derivations)  by the even operator
\[
E_\alpha:=t\partial_t+\alpha\theta_1\partial_{\theta_1}+\theta_2\partial_{\theta_2}.
\]

There is only one non-trivial central extension (for a proof, see \cite{KvL}) given by the even cocycle $c$ with the following non-zero values:
\[
\begin{array}{ll}
c(L_m, L_n)&=\nfrac12m(m^2-(\alpha+1)^2)\delta_{m,-n},\\

c(t^{n+1}\partial_{\theta_i}, S_m^j)&=-m(m-\alpha+1)\delta_{m,-n}\delta_{m,j},\\

c(G_m, G_n)&= m\delta_{m,-n},\\

c(t^{m}\theta_1\partial_{\theta_2}, t^{n}\theta_2\partial_{\theta_1})&=m\delta_{m,-n},\\
\end{array}
\]
where
\[
S_n^j:=t^{n-1}\theta_j\left(t\partial_t +\nfrac12(n+\alpha)(\theta_1\partial_{\theta_1}+\theta_2\partial_{\theta_2})\right).
\]

\sssbegin{Theorem}[No NIS on $\widehat{\fsvect}_\alpha^L(1|2)$]\label{noNIS}
On $\widehat{\fsvect}_\alpha^L(1|2)$ with any of Cartan matrices \eqref{eq-3d-137}, there is no NIS compatible with the principal grading $\deg X_i^\pm=\pm1$ for all $i$. \emph{Conjecturally, there is no NIS whatsoever}.\end{Theorem}

\begin{proof} Induction, beginning with degree 0, see \eqref{recipe}, leads to a contradiction. \end{proof}

\ssec{The exceptional simple vectorial Lie superalgebras for $p=3$}\label{SSp=3}

\sssec{$\fv\fas(4;\protect\un|4)$ for $p=3$}\label{sssVAS} For $\fg=\fv\fas(4;\un|4)$ and any $p\neq 2$, we have the following description:
$\fg_\ev=\fvect(4;\un|0)$,
and $\fg_\od=\Omega^{1}(4;\un|0)\otimes_{\Omega^0(4;\un|0)} \Vol^{-1/2}(4;\un|0)$ with the
natural $\fg_\ev$-action on $\fg_\od$ and the bracket of odd
elements is given by
\[
\left[\frac{\omega_{1}}{\sqrt{\vvol}},
\frac{\omega_{2}}{\sqrt{\vvol}}\right]=
\frac{d\omega_{1}\wedge\omega_{2}+ \omega_{1}\wedge d\omega_{2}}{
\vvol},
\]
where we identify
\[
\frac{dx_{i}\wedge dx_{j}\wedge dx_{k}}{
\vvol}=\sign(ijkl)\partial_{x_{l}}\text{ for any permutation $(ijkl)$
of $(1234)$}.
\]
For $p=3$, when $-\nfrac12=1$, we have the following odd NIS on $\fv\fas(4;\un|4)$:
\[
(f(u)\partial_i,\ g(u)du_j\vvol)=\delta_{i,j}\int fg \vvol.
\]

Having in mind  Fact (Subsection~\ref{Fact}) and the fact that $\fv\fas$ is a deform of $\fs\fle(4)$, it is no wonder that there is an odd NIS on $\fv\fas$.

\sssec{Lie superalgebras indigenous to $p=3$}\label{SSlopsup=3} In~\cite{BoL1,BoGL}, there are constructed two super analogs of the Melikyan algebra (one, as a Cartan prolong of a non-positive part of the exceptional Lie superalgebra $\fg(2,3)$, the other one analogous of the original Melikyan's construction, with $\fk(1;\One|1)$ for $p=3$ instead of $\fvect(1;\One)$ for $p=5$), and 5 more vectorial algebras indigenous to $p=3$. Lemma~\ref{Pr 4.7} shows that the only superalgebras among the following seven exceptions
\footnotesize
\begin{equation}\label{BJ}
\begin{tabular}{|l|l|}
\hline
$\mathfrak{Bj}(1;\un|7)$&  see \cite{BoL1}\\
\hline
$\mathfrak{Bj}(3;\un|5)$& $\un$ was forgotten in \cite{BoGL}, recovered in \cite{BGLLS1}\\
\hline
$\mathfrak{Me}(3;\un|3)$,\ $\mathfrak{Bj}(3;\un|3)$& see \cite{BoGL}\\
\hline
$\mathfrak{Bj}(4;\un|5)$& missed in \cite{BoGL}, added in \cite{BGLLS1} \\
\hline
$\mathfrak{Bj}(3;\un|4)$& see \cite{BGLLS1} (denoted by $\mathfrak{BRJ}$ in \cite{BoGL})  \\
\hline
$\mathfrak{Brj}(4|3)$& see the latest arXiv version of \cite{BoGL} \\
\hline
\end{tabular}
\end{equation}
\normalsize
that could have a NIS are $\mathfrak{Me}(3;\un|3)$, $\mathfrak{Bj}^{(2)}(3;\un|5)$, $\mathfrak{Bj}^{(1)}(3;\un|4)$, $\mathfrak{Brj}^{(1)}(4|3)$, $\mathfrak{Bj}(1;\un|7)$, and $\fbj=\fg(2,3)^{(1)}/\fc$. The last one has a NIS thanks to~\eqref{recipe}, whereas \textit{SuperLie} proves that each of the other algebras has a NIS. These NISes are corollaries of the \textit{necessary}, but not sufficient, conditions of Lemma~\ref{Pr 4.7};  the ambient algebras do not have NIS. We give analogs of eqs.~\eqref{bilFvect1} and~\eqref{bilFvect2}, which also have to be proved separately, in Subsection~\ref{ssAbo}.
\textbf{For simplicity, we give the answer only for $\un=\One$}.

$\bullet$ Observe that $\mathfrak{Me}(3;\un|3)$ is a subalgebra of $\fk(t, x_1, x_2;\un|x_3, x_4, x_5)$, the latter
preserving to the distribution singled out by the form (see~\cite{BoGL})
\[
\alpha= dt- x_1dx_2 + x_2dx_1 + x_3dx_3 + x_4dx_4 + x_5dx_5.
\]
On $\fk(3;\un|3)$, there is no NIS, see Subsection~\ref{sssK}. The odd NIS on $\fg=\mathfrak{Me}(3;\One|3)\subset \fk(3;\un|3)$ is given by pairing $\fg_{-2}$ with the component $\fg_h$ of highest degree $h$, and extending to other pairs of elements by invariance. For example, for $\un=\One$ we have $h=5$ and
\[
(1,\ 2t^{(2)} x_2^{(2)} x_3 + t^{(2)} x_1 x_2 x_4 + 2 t^{(2)} x_1^{(2)} x_5
              + 2 t^{(2)} x_3 x_4 x_5 + x_1^{(2)} x_2^{(2)} x_3 x_4 x_5) = 1.
\]

$\bullet$  For any $p\neq 2$ and Cartan matrix of $\fag(2)$, the Cartan prolongs of the non-positive part in any its $\Zee$-grading return
$\fag(2)$, except for $p=3$ where for the Cartan matrix~$1)$  in~\eqref{ag2}
and its $\Zee$-grading $r=(100)$, the Cartan prolongs yield the simple
subalgebras of series $\mathfrak{Bj}(1;\un |7)$, and the
exceptional algebra  $\mathfrak{bj}=\fg(2,3)^{(1)}/\fc$ whose double extension $\fg(2,3)$ has 5 inequivalent Cartan matrices, see \cite{BGL}.

The Lie superalgebra $\fg=\fag(2)$ over $\Cee$ has the
following four non-equivalent Cartan matrices, see~\cite{CCLL, BGL}; we use Chevalley generators corresponding to matrix 1):
\tiny \be\label{ag2}
1)\; \begin{pmatrix} 0 & 1 & 0 \\ -1 & 2 & -3 \\ 0 & -1 & 2
\end{pmatrix}\quad 2)\; \begin{pmatrix}
0 & 1 & 0 \\ -1 & 0 & 3 \\ 0 & -1 & 2
\end{pmatrix}\quad 3)\; \begin{pmatrix}
0 & -3 & 1 \\ -3 & 0 & 2 \\ -1 & -2 & 2
\end{pmatrix}\quad 4) \; \begin{pmatrix}
2 & -1 & 0 \\ -3 & 0 & 2 \\ 0 & -1 & 1
\end{pmatrix}
\ee
\normalsize
From an explicit description of $\fg(2)$ and its first fundamental
representation, e.g., see~\cite{FH}, we deduce an explicit form of the
non-positive elements of $\fg\subset\fk(1|7)$ which we give in terms of the
generating functions with respect to the contact bracket
corresponding to the contact form
\begin{equation}\label{alpha}
\alpha=dt-\sum\limits_{i=1,2,
3}(v_idw_i+w_idv_i)+2udu,
\end{equation}
where the $v_i$, $w_i$ and $u$
are odd; our notation match~\cite{FH}, p. 354, while our $X_i^+$
and $X_i^-$ correspond to $X_i$ and $Y_i$ of~\cite{FH}, p.~340, respectively.

We also set $X_3^\pm :=[X_1^\pm , X_2^\pm ]$, $X_4^\pm
:=[X_1^\pm , X_3^\pm ]$, $X_5^\pm :=[X_1^\pm , X_4^\pm ]$, $X_6^\pm :=[X_2^\pm , X_5^\pm ]$.

To describe the $\fg_0$-module $\fg_{-1}=\Span(u\;\text{ and $v_i,
w_i$ for $i=1,2, 3$})$, only the highest weight vector suffices:
\[\footnotesize
\renewcommand{\arraystretch}{1.4}
\begin{tabular}{|l|l|}
\hline
$\fg_{i}$&the generating functions of generators of the $\fg_{0}$-modules $\fg_i$\\
\hline \hline
$\fg_{-2}$&$1$\\
\hline
$\fg_{-1}$& $v_3$\\
\hline
$\fg_0$&$\fz=\Span(t)$,\ \
$X_1^+= -v_3w_2-uv_1,\quad X_1^-=-v_2w_3-uw_1$,\ \
$X_2^+=v_2w_1,\quad X_2^-=v_1w_2$.\\
\hline
\end{tabular}
\]
As expected, for $p=0$ and $p>3$, the CTS prolong is isomorphic to $\fag(2)$.

For $p=3$, the Lie algebra $\fg_0$ is not simple, but has a simple
Lie subalgebra isomorphic to $\fpsl(3)$ generated by
$x_1^\pm=X_1^\pm$, and $x_2^\pm =[X_{1}, X_2^\pm]$.
Let
$\widetilde\fg_0:=\fpsl(3)\oplus\fz$, where the center of
$\widetilde\fg_0$ is $\fz=\Span(t +
v_{1}w_{1} +  v_{2} w_{2} + 2\,  v_{3} w_{3})$.

The $\widetilde\fg_0$-module $\fg_{1}$ splits into two irreducible
components: a $0|1$-dimensional, and $0|7$-dimensi\-o\-nal with lowest
weight vectors, respectively:
\[\footnotesize
\renewcommand{\arraystretch}{1.4}
\begin{tabular}{|l|l|}
\hline
$V_{1}'$&$v_1v_2w_3+v_1uw_1+v_2uw_2+2v_3uw_3+v_3w_1w_2$\\
$V_{1}''$& $tw_3+v_1w_1w_3+v_2w_2w_3+uw_1w_2$\\
\hline
\end{tabular}
\]

Since $\fg_{1}$ generates the positive part of the CTS prolong,
$[\fg_{1}, \fg_{-1}]=\widetilde\fg_{0}$, the
$\widetilde\fg_0$-module $\fg_{-1}$ is irreducible, and
$[\fg_{-1}, \fg_{-1}]=\fg_{-2}$, the standard criterion for
simplicity
ensures that the generalized Cartan prolong is simple; we denote it
$\mathfrak{Bj}(1;\un|7)$.

The positive components of $\mathfrak{Bj}(1;\un|7)$ are all of dimension
8 (the direct sums $\fg_k=\fg_k'\oplus\fg_k''$ of irreducible
$\fg_0$-modules of dimension 1 and 7, up to parity), except the one of
the highest degree ($2(3^N-1)+3-2=2\cdot 3^N-1$), which is of dimension 1,
and the second highest degree, which is of dimension 7. Let $V_{k}'$ and
$V_{k}''$ be the lowest weight vectors (with respect to $\fg_0$) of $\fg_k'$ and $\fg_k''$,
respectively for $\un=1$:\footnotesize
\[
\renewcommand{\arraystretch}{1.4}
\begin{tabular}{|l|l|}
\hline
$V_{2}'$&$t^2+2v_3uw_1w_2+v_2v_3w_2w_3+v_1v_3w_1w_3+2v_1v_2w_1w_2+2v_1v_2uw_3$\\
$V_{2}''$& $2v_1uw_1w_3+2v_2uw_2w_3+2v_3w_1w_2w_3+tuw_3+2tw_2w_3$\\
\hline
$V_{3}'$&$tv_1v_2w_3+tv_1uw_1+tv_2uw_2+2tv_3uw_3+tv_3w_1w_2$\\
$V_{3}''$&
$t^2w_3+tuw_1w_2+2v_1v_2w_1w_2w_3+tv_1w_2w_3+tv_2w_2w_3+2v_3uw_2w_3$\\
\hline
$V_{4}''$&
$2v_1v_2uw_1w_2w_3+2tv_1uw_1w_3+2tv_2uw_2w_3+2tv_3w_1w_2w_3+t^2uw_3+2t^2w_1w_2$\\
\hline
$V_{5}'$&$2v_1v_2v_3uw_1w_2w_3+2t^2v_1v_2w_3+2t^2v_1uw_1+
2t^2v_2uw_2+t^2v_3uw_3+2t^2v_3w_1w_2
$\\
\hline
\end{tabular}
\]\normalsize

A subalgebra of $\mathfrak{Bj}(1;\un|7)$: let $\fg_1''$ be generated by  $V_1''$, as
$\widetilde\fg_0$-module. Let $\mathfrak{bj}:=(\fg_-, \fg_0, \fg_1'')_*$ denote the
partial prolong. This is a simple Lie superalgebra  of dimension
$(10|14)$ with
$$
\mathfrak{bj}_2=\Span(t^2+2v_1v_2uw_2+2v_1v_2w_1w_2
+v_1v_2w_1w_3+v_2v_3w_2w_3+2v_3uw_1w_2)
$$
and $\mathfrak{bj}_3=0$.

\ssec{Abomination expressions}\label{ssAbo} In what follows  elements of different parity are separated by a bar: (even$|$odd).

$\bullet$ The odd NIS on $\mathfrak{Bj}(1;\One|7)$ is given by paring of $\fg_{-2}$ with $\fg_{5}$, and extending to other pairs of elements by invariance:
\[
  (-1,\ \eta _1 \eta _2 \eta _3 \theta _1 \xi _1 \xi _2 \xi _3+2 \eta _1 \theta _1 \xi _1 t^{(2)}+\eta _2 \theta _1
  \xi _2 t^{(2)}+\eta _3 \theta _1 \xi _3 t^{(2)}+\eta _2 \eta _3 \xi _1 t^{(2)}+\eta _1 \xi _2 \xi _3
  t^{(2)})  =1.
\]

$\bullet$ We realize~$\mathfrak{Brj}^{(1)}(3;\One;|4)$ as a Lie subsuperalgebra of~$\fvect(3;\One|4)$ in
indeterminates  (even$|$odd) $x_1,x_2,x_3\mid x_4,x_5,x_6,x_7$.
The odd NIS on $\mathfrak{Bj}^{(1)}(3;\One|4)$ is given by paring of $\fg_{-3}$ with
$\fg_{6}$, and extending to other pairs of elements by invariance, where  $p(v) = \od$  and $\fg_{-3}=\Span(v = \del_4)$; for a basis in~$\fg_6$ we take the even element $w$:
\tiny
\[
\begin{array}{ll}
(v,\ w)& = 1\\
  w =
  {}&x_1 x_2^{(2)} x_3\del_1+x_1^{(2)} x_3^{(2)}\del_1+2 x_1 x_2 x_3^{(2)}\del_2+2 x_2^{(2)} x_3^{(2)}\del_3
      +x_2^{(2)} x_4\del_5+2 x_1 x_3 x_4\del_5+2 x_2 x_3 x_4\del_6+2 x_2^{(2)} x_3 x_4\del_4+2 x_3^{(2)} x_4\del_7\\
  {}&{}+x_1 x_3^{(2)} x_4\del_4+x_1 x_3^{(2)} x_5\del_5+2 x_2 x_3^{(2)} x_5\del_6+2 x_3^{(2)} x_4 x_5\del_1
      +2 x_1 x_2 x_3 x_6\del_5+2 x_2^{(2)} x_3 x_6\del_6+2 x_1 x_3^{(2)} x_6\del_6+2 x_2 x_3^{(2)} x_6\del_7\\
  {}&{}+2 x_2 x_3 x_4 x_6\del_1+2 x_3^{(2)} x_4 x_6\del_2+2 x_2 x_3^{(2)} x_5 x_6\del_1+2x_1 x_2^{(2)}x_7\del_5
      +2 x_1^{(2)} x_3 x_7\del_5+x_1 x_2 x_3 x_7\del_6+x_2^{(2)} x_3 x_7\del_7+x_3^{(2)} x_4 x_7\del_3\\
  {}&{}+x_3 x_4 x_5 x_7\del_5+x_3^{(2)} x_4 x_5 x_7\del_4+x_1 x_2 x_3 x_6 x_7\del_1
      +x_2 x_3^{(2)} x_6 x_7\del_3+2 x_2 x_4 x_6 x_7\del_5+x_3 x_4 x_6 x_7\del_6
      +2 x_3^{(2)} x_5 x_6 x_7\del_6\\
  {}&{}+2 x_2 x_3^{(2)} x_5 x_6 x_7\del_4.
\end{array}
\]
\normalsize

$\bullet$ We realize~$\mathfrak{Brj}^{(2)}(3;\One;|5)$ as a Lie subsuperalgebra of~$\fvect(3;\One|5)$ in
indeterminates (even$|$odd) $x_1,x_4,x_5\mid x_2,x_3,x_6, x_7,x_8$. The  odd NIS on
$\mathfrak{Brj}^{(2)}(3;\One|5)$ is given by paring of $\fg_{-2}$ with $\fg_{5}$,
and extending to other pairs of elements by invariance, where $\fg_5=\Span(w_2, w_3\mid w_1)$ and $\fg_{-2}=\Span(v_1 = \del_1\mid v_2 = \del_2,\ \ v_3 = \del_3)$:
\tiny{
\[
\begin{array}{ll}
(v_1,\ w_1) &= 2,\quad (v_2,\ w_3) = 1,\quad (v_3,\ w_2) = 2,\\

w_1 =
  {}& (x_1^{(2)} x_2)\del_4+2 (x_1^{(2)} x_3)\del_5 +2 (x_1^{(2)} x_6  x_7)\del_6
      +2 (x_1 x_3 x_6 x_7)\del_4+(x_1^{(2)} x_6 x_8)\del_7+2(x_1 x_3 x_6 x_8)\del_5\\
      
  {}&{}+2 (x_1^{(2)} x_7 x_8)\del_8+2 (x_1 x_2 x_7 x_8)\del_5
      +2 (x_1^{(2)} x_6 x_7 x_8)\del_1+2 (x_1 x_2 x_6 x_7  x_8)\del_2
      +2 (x_1 x_3 x_6 x_7 x_8)\del_3,\\
      
  w_2 =
  {}& 2 (x_1^{(2)} x_2)\del_6+2 (x_1^{(2)} x_3)\del_7+2 (x_1 x_2  x_3)\del_4
      +(x_1^{(2)} x_3 x_4)\del_2+2 (x_1^{(2)} x_2 x_5)\del_2+(x_1^{(2)} x_3 x_5)\del_3+2 (x_1^{(2)} x_4  x_5)\del_4+(x_1^{(2)} x_5^{(2)})\del_5\\
      
  {}&{}+(x_1^{(2)} x_5 x_6)\del_6
      +2(x_1 x_3 x_5 x_6)\del_4+(x_1^{(2)} x_3 x_7)\del_1+(x_1^{(2)} x_4  x_7)\del_6+2 (x_1 x_3 x_4 x_7)\del_4+(x_1 x_2 x_5 x_7)\del_4+2(x_1 x_3 x_5 x_7)\del_5\\
      
  {}&{}
      +(x_1 x_3 x_5 x_6 x_7)\del_2+(x_1 x_5^{(2)} x_6  x_7)\del_4+2 (x_1^{(2)} x_4 x_8)\del_7+2 (x_1 x_3 x_4  x_8)\del_5+2 (x_1^{(2)} x_5 x_8)\del_8
      +(x_1 x_2 x_5  x_8)\del_5\\
      
  {}&{}+(x_1 x_3 x_4 x_6 x_8)\del_2+2 (x_1 x_2 x_5 x_6  x_8)\del_2+(x_1 x_5^{(2)} x_6 x_8)\del_5+2 (x_1 x_2 x_7  x_8)\del_7
      +2 (x_1 x_3 x_7 x_8)\del_8+(x_2 x_3 x_7  x_8)\del_5\\
      
  {}&{}+(x_1^{(2)} x_4 x_7 x_8)\del_1+2 (x_1 x_3 x_4 x_7 x_8)\del_3+(x_1 x_2 x_5 x_7 x_8)\del_3
      +(x_1 x_4 x_5 x_7  x_8)\del_5+(x_2 x_3 x_6 x_7 x_8)\del_2+(x_1 x_5 x_6 x_7 x_8)\del_7\\
      
  {}&{}+(x_3 x_5 x_6 x_7 x_8)\del_5+(x_1 x_4 x_5 x_6 x_7  x_8)\del_2
      +(x_1 x_5^{(2)} x_6 x_7 x_8)\del_3,\\
      
  w_3 = {}& 2 (x_1^{(2)} x_2)\del_7+2 (x_1^{(2)} x_3)\del_8+(x_1 x_2 x_3)\del_5+(x_1^{(2)} x_2 x_4)\del_2
            +2 (x_1^{(2)} x_3 x_4)\del_3+(x_1^{(2)} x_4^{(2)})\del_4+(x_1^{(2)} x_2 x_5)\del_3+2 (x_1^{(2)} x_4 x_5)\del_5\\
            
  {}&{}+2 (x_1^{(2)} x_3 x_6)\del_1
      +2 (x_1^{(2)} x_4 x_6)\del_6+(x_1 x_3 x_4 x_6)\del_4+2 (x_1^{(2)} x_5 x_6)\del_7+2 (x_1 x_2 x_5 x_6)\del_4+(x_1^{(2)} x_2 x_7)\del_1
      +2 (x_1 x_2 x_4 x_7)\del_4\\
  {}&{}+(x_1 x_3 x_4 x_7)\del_5+(x_1^{(2)} x_5 x_7)\del_8+2 (x_1 x_2 x_5 x_7)\del_5+(x_1 x_2 x_6 x_7)\del_6
      +(x_1 x_3 x_6 x_7)\del_7+(x_2 x_3 x_6 x_7)\del_4+2 (x_1 x_3 x_4 x_6 x_7)\del_2\\
  {}&{}+(x_1^{(2)} x_5 x_6 x_7)\del_1+(x_1 x_2 x_5 x_6 x_7)\del_2
      +2 (x_1 x_4 x_5 x_6 x_7)\del_4+(x_1^{(2)} x_4 x_8)\del_8+2 (x_1 x_2 x_4 x_8)\del_5+2 (x_1 x_2 x_6 x_8)\del_7\\
  {}&{}+2 (x_1 x_3 x_6 x_8)\del_8
      +(x_2 x_3 x_6 x_8)\del_5+(x_1 x_2 x_4 x_6 x_8)\del_2+2 (x_1 x_4 x_5 x_6 x_8)\del_5+2 (x_1 x_2 x_4 x_7 x_8)\del_3+2 (x_1 x_4^{(2)} x_7 x_8)\del_5\\
  {}&{}
      +(x_1 x_2 x_6 x_7 x_8)\del_1+(x_2 x_3 x_6 x_7 x_8)\del_3+2 (x_1 x_4 x_6 x_7 x_8)\del_7+2 (x_3 x_4 x_6 x_7 x_8)\del_5
      +2 (x_1 x_4^{(2)} x_6 x_7 x_8)\del_2\\
  {}&{}+2 (x_1 x_5 x_6 x_7 x_8)\del_8+(x_2 x_5 x_6 x_7 x_8)\del_5+2 (x_1 x_4 x_5 x_6 x_7 x_8)\del_3.
\end{array}
\]
}\normalsize

$\bullet$ We realize~$\mathfrak{Brj}^{(1)}(4|3)$ as a Lie subsuperalgebra of~$\fvect(4;\One|3)$ in
indeterminates (even$|$odd) $x_1,x_2,x_3,x_4\mid x_5,x_6,x_7$.
The odd NIS on $\mathfrak{Brj}^{(1)}(4|3)$ is given by paring of $\fg_{-4}$ with $\fg_{10}$, and extending to other pairs of elements by invariance,
where $\fg_{-4}=\Span(v_2 = 2\del_1\mid v_1 = \del_5)$ and $\fg_{10}=\Span(w_2\mid w_1)$:
{\tiny
\[
\begin{array}{ll}
 (v_1,\ w_2) &= 1,\quad (v_2, w_1) = 2,\\
   w_1 =
  {}&2 x_1^{(2)} x_2\del_4+x_1^{(2)} x_2^{(2)}\del_1+x_1^{(2)} x_2 x_3\del_2+2 x_1 x_2^{(2)} x_3\del_3
      +2 x_1^{(2)} x_3^{(2)}\del_3+x_1 x_2^{(2)} x_4\del_4+2 x_1^{(2)} x_3 x_4\del_4+2 x_1^{(2)} x_5\del_6+x_1x_2 x_5\del_7\\
  {}&{}+x_1 x_2^{(2)} x_5\del_5+2 x_1^{(2)} x_3 x_5\del_5+2 x_1 x_2 x_4 x_5\del_6+2 x_2^{(2)} x_4 x_5\del_7
      +x_1 x_3 x_4 x_5\del_7+x_2^{(2)} x_4^{(2)} x_5\del_6+2 x_1 x_3 x_4^{(2)} x_5\del_6\\
  {}&{}+2 x_1 x_2^{(2)} x_6\del_6+x_1 x_2 x_3 x_6\del_7+2 x_1 x_2 x_3 x_4 x_6\del_6+2 x_2^{(2)}x_3x_4x_6\del_7
      +2 x_1 x_3^{(2)} x_4 x_6\del_7+x_2^{(2)} x_3 x_4^{(2)} x_6\del_6+x_1 x_3^{(2)} x_4^{(2)} x_6\del_6\\
  {}&{}+2 x_1^{(2)} x_5 x_6\del_1+x_1 x_2 x_5 x_6\del_2+2 x_2^{(2)} x_5 x_6\del_3+x_1 x_3 x_5 x_6\del_3
      +x_1 x_2 x_4 x_5 x_6\del_1+2 x_2^{(2)} x_4^{(2)} x_5 x_6\del_1+x_1 x_3 x_4^{(2)} x_5 x_6\del_1\\
  {}&{}+x_1 x_2^{(2)} x_7\del_7+x_1^{(2)} x_3 x_7\del_7+2 x_1 x_2 x_3 x_4 x_7\del_7+x_2^{(2)} x_3 x_4^{(2)}x_7\del_7
      +x_1 x_3^{(2)} x_4^{(2)} x_7\del_7+2 x_1 x_2 x_5 x_7\del_3+x_1 x_2 x_4 x_5 x_7\del_2\\
  {}&{}+x_2^{(2)} x_4 x_5 x_7\del_3+2 x_1 x_3 x_4 x_5 x_7\del_3+2 x_1 x_2 x_4^{(2)} x_5 x_7\del_1
      +x_2^{(2)} x_4^{(2)} x_5 x_7\del_2+2 x_1 x_3 x_4^{(2)} x_5 x_7\del_2+x_1^{(2)} x_6 x_7\del_4\\
  {}&{}+2 x_1 x_2 x_3 x_6 x_7\del_3+x_1 x_2 x_4 x_6 x_7\del_4+x_1 x_2 x_3 x_4 x_6 x_7\del_2
      +x_2^{(2)} x_3 x_4 x_6 x_7\del_3+x_1 x_3^{(2)} x_4 x_6 x_7\del_3+2 x_2^{(2)} x_4^{(2)} x_6 x_7\del_4\\
  {}&{}+x_1 x_3 x_4^{(2)} x_6 x_7\del_4+2 x_1 x_2 x_3 x_4^{(2)} x_6 x_7\del_1+x_2^{(2)} x_3 x_4^{(2)} x_6 x_7\del_2
      +x_1 x_3^{(2)} x_4^{(2)} x_6 x_7\del_2+2 x_1 x_5 x_6 x_7\del_7+x_1 x_2 x_5 x_6 x_7\del_5\\
  {}&{}+2 x_2 x_4 x_5 x_6 x_7\del_7+2 x_2^{(2)} x_4 x_5 x_6 x_7\del_5+x_1 x_3 x_4 x_5 x_6 x_7\del_5+2 x_3 x_4^{(2)} x_5 x_6 x_7\del_7,\\
  w_2 =
  {}&x_1^{(2)} x_2\del_7+x_1^{(2)} x_2^{(2)}\del_5+2 x_1^{(2)} x_2 x_4\del_6+2 x_1 x_2^{(2)} x_4\del_7
      +x_1 x_2^{(2)} x_4^{(2)}\del_6+2 x_1^{(2)} x_5\del_3+x_1^{(2)} x_2 x_6\del_2+2 x_1 x_2^{(2)} x_6\del_3\\
  {}&{}+2 x_1^{(2)} x_4 x_6\del_4+x_1^{(2)} x_2 x_4 x_6\del_1+2 x_1 x_2^{(2)} x_4^{(2)} x_6\del_1
      +2 x_1 x_4 x_5 x_6\del_7+x_1 x_4^{(2)} x_5 x_6\del_6+2 x_1^{(2)} x_2 x_7\del_3+x_1^{(2)} x_2 x_4x_7\del_2\\
  {}&{}+x_1 x_2^{(2)} x_4 x_7\del_3+x_1^{(2)} x_4^{(2)} x_7\del_4+2 x_1^{(2)} x_2 x_4^{(2)} x_7\del_1
      +x_1 x_2^{(2)} x_4^{(2)} x_7\del_2+x_1 x_4^{(2)} x_5 x_7\del_7+x_1^{(2)} x_6x_7\del_7+x_1^{(2)}x_2x_6x_7\del_5\\
  {}&{}+2 x_1 x_2 x_4 x_6 x_7\del_7+2 x_1 x_2^{(2)} x_4 x_6 x_7\del_5+x_2^{(2)} x_4^{(2)} x_6 x_7\del_7
      +x_1 x_4 x_5 x_6 x_7\del_3+x_1 x_4^{(2)} x_5 x_6 x_7\del_2.
\end{array}
\]
}

%

\section{Summary: NIS on simple Lie (super)algebras}

 Some NISes listed below under even (resp. odd) Subsection are actually odd (resp. even) depending on a parameter, as indicated.

\ssec{NIS even}  1) Lie (super)algebras $\fg(A)$ over $\Cee$ with any symmetrizable Cartan matrix $A$ have a NIS; these algebras are simple if and only if $A$ is invertible. (Aside: These algebras are $\Zee$-graded, and as such they fall into three classes: finite-dimensional, of polynomial growth and of exponential growth. The simple Lie superalgebras of the first two types, and certain \lq\lq islands" in the ocean of the third type are classified, see \cite{CCLL} and references therein.)

2) $\dim\fg<\infty$ over $\Cee$. There is an even NIS on $\fh^{(1)}(2n)$. Its double extension is $\fpo(0|2n)$.

2p) Over $\Kee$ of characteristic $p$, there is a NIS on $\fh^{(1)}(0|m)$ for $m\geq 4$, and  we see that $p(\text{NIS})=m\pmod 2$. The double extension of $\fh^{(1)}(0|m)$ is $\fpo(0|m)$ except for $p=2$ and $m=4$, where there are 3 non-isometric double extensions, see \cite{BeBou}.

3) (Twisted) loop superalgebras over $\Cee$.  These are  Lie superalgebras $(\fg)_\psi^{\ell(m)}$ corresponding to lines 1 through 10  of table~\eqref{**}. (Aside: The twisted loop superalgebras corresponding to lines 5, 6, 8, 9 of table~\eqref{**} and loop superalgebras $(\fh^{(1)}(2n))^{\ell(1)}$  do not have any Cartan matrix.)

4) Stringy Lie superalgebra $\fk^L(1|6)$ over $\Cee$ has a NIS. No central extensions, no outer derivations.

5) The following simple vectorial Lie superalgebras in positive characteristic have NIS:

$\fk(2n+1;\un|m)$ if $2n+2-m\equiv -4\pmod p$, see \eqref{cond},
where $p>2$, and $p(\text{NIS})=m\pmod2$;

the Melikyan algebra (for $p=5$);

several of Skryabin algebras for $p=3$, namely, $\fby^{(1)}(7;\un)$, $\fmy(6; \un)$, $\fs\fmy(6; \un)$, see \cite{GL};

the Lie algebra $\fvect(n;\un)$ has a NIS if and only if either $n=1$ and $p=3$
or $n=p=2$; if $n>2$, then $\fsvect^{(1)}(n;\un)$ has a NIS if and only if $n=3$, see Subsection~\ref{Th6.4};

the simple Lie algebra~$\fh^{(\infty)}_\omega(2k;\un|0)$ has a NIS, see Subsection~\ref{hNIS};
 
for deforms of $\fo(5)$ for $p=3$, see Claim~\ref{cl_o(5)};

for deforms of Lie algebras $\fg=\fbr(3)$ and Lie superalgebra $\fbrj(2;3)$ for $p=3$, and $\fwk(4;\alpha)$, where $\alpha\neq 0, 1$, for $p=2$, see Claim~\ref{wk4};

for deforms of Lie algebras $\fwk(3;\alpha)$ and $\fwk^{(1)}(3;\alpha)/\fc$, where $\alpha\neq 0, 1$, for $p=2$, see Claim~\ref{wk3}.

\ssec{NIS odd}  1) $\dim\fg<\infty$ over $\Cee$. These are $\fpsq(n)$ for $n\geq 3$ and
$\fh^{(1)}(0|2n+1)$ for $n\geq 2$. The corresponding double extensions are $\fq(n)$ and
$\fpo(0|2n+1)$.

1p) Over $\Kee$ of characteristic $2$ and $m\geq 5$ odd, where there are 2 non-isometric double extensions of $\fh^{(1)}(0|m)$, see \cite{BeBou}.

2) (Twisted) loop superalgebras over $\Cee$. These are  $(\fpsq(n))^{\ell(2)}$ and $(\fh^{(1)}(0|2n+1))^{\ell(2)}$ and Lie superalgebras corresponding to lines 10 and 11 of table~\eqref{**}. Observe that a (there are several, their isometry classes are unknown) double extension corresponding to line 10 does have Cartan matrices, albeit non-symmetrizable ones.

3) Stringy Lie superalgebra $\fk^L(1|7)$  over $\Cee$ has a NIS. No central extensions, no outer derivations.

Observe that $\widehat\fsvect_\alpha^L(1|2)$ has a Cartan matrix, albeit non-symmetrizable one, and $\fsvect_\alpha^L(1|2)$ has a non-trivial central extension and an outer derivation, but no NIS; no NIS on $\widehat\fsvect_\alpha^L(1|2)$, either.

4)  For $p=2$, the queerification $\fq(\fg)$ (see \cite{BLLSq}) of any simple restricted Lie algebra $\fg$ with a NIS has an odd NIS, see Subsection  \ref{sssQg}.

5) There is a NIS on $\mathfrak{sb}(n;\un)$ and $\mathfrak{sle}^{(1)}(n;\un)$, where $p(\text{NIS})\equiv n+1\pmod 2$, see Subsection~\ref{ssSB}.

6) On $\mathfrak{Me}(3;\un|3)$, $\mathfrak{Bj}^{(2)}(3;\un|5)$, $\mathfrak{Bj}^{(1)}(3;\un|4)$, $\mathfrak{Brj}^{(1)}(4|3)$, $\mathfrak{Bj}(1;\un|7)$, see Subsection~\ref{SSlopsup=3}.

7) on $\fv\fas(4;\un|4)$ for $p=3$, see Subsection~\ref{sssVAS}.

\ssec{Open problems}\label{OpP} 1) Extend  to any $\un$ the results and formulas obtained in Subsections~\ref{SSlopsup=3} and \ref{ssAbo} for $\un=\One$.

2) Conjecturally, all non-trivial central extensions of (twisted) loop superalgebras with values in finite-dimensional simple Lie superalgebras are those listed in Subserction~\ref{CentExtLoo}. Which of these extensions can be extended to double extensions and which of these double extensions belong to one isometry class?

3) Prove conjectures \ref{ClSvect}, \ref{conj1} and \ref{noNIS}.

4) Investigate if there is a NIS on deforms of $\fb_{\lambda}(n)$, see Subsection~\ref{thdefb}.

5) Related to this paper are several conjectures, listed in the order of feasibility. 

5a) Classification of simple $\Zee$-graded Lie superalgebras of polynomial growth over $\Cee$: they are (a) the finite-dimensional ones, (b) vectorial (with polynomial coefficients), (c) stringy (with Laurent coefficients) and (d) (twisted) loop superalgebras. Since the time this conjecture was formulated in \cite{LSS} it was partly proved. Still an \textbf{open problem} is to prove that the  (twisted) loop Lie superalgebras~\eqref{**} are all that remains to add to the classifications in cases (a)--(c).

5b) Classification of simple finite-dimensional Lie algebras over an algebraically closed field $\Kee$ of characteristic $3$: \textbf{conjecturally} these are the examples obtained by means of the Kostrikin-Shafarevich procedure, examples listed in \cite{GL}, and deforms of these two types of examples.

5c) Classification of simple finite-dimensional Lie superalgebras over an algebraically closed field $\Kee$ of characteristic $>5$. For a formulation of the conjecture, see \cite{BGLLS}.

6) Consider deforms with odd parameters; for the cases where they are classified, see \cite{BGLd}.

7) Consider general Lie algebras of matrices of complex size $\fgl(\lambda):=(U(\fsl(2))/(C_2-\lambda^2+1))_L$, where $A_L$ is the Lie algebra associated with the associative algebra $A$ and $C_2$ is the 2nd order Casimir, and the ortho/sym\-plectic subalgebras of $\fgl(\lambda)$, where $\lambda \in\Cee$, see \cite{GL96, LSe}. Consider generalizations of these Lie algebras related with simple Lie algebras $\fg$ of rank $>1$ instead of $\fsl(2)$.
Consider super versions of these Lie algebras. Most of them (or all?) have a NIS.

8) Simple Lie (super)algebras for $p= 2$. For a general theory of double extensions ---  quite different from that for $p\neq 2$ --- and interesting examples, see  \cite{BeBou}. Additionally, for selected examples, see Subsections~\ref{Asym}, \ref{wk4}, \ref{wk3}; for a wide series of examples, see Subsection~\ref{sssQg}.

9) Lie algebras of pseudodifferential operators, see \cite{Dz4}.

\vfill
\newpage

\section{Tables}

In Tables~\eqref{*} and~\eqref{**}, let $\Ad_A(X):=A^{-1}XA$ for any $X\in\fgl(m|n)$ and some even invertible $A\in\fgl(m|n)$. 
\footnotesize

\begin{table}[ht]\centering
{%
\caption{Finite order automorphisms of simple finite-dimensional Lie superalgebras $\fg$ over $\Cee$ and of their double extensions (\cite{Se})}\label{td4.0}\nopagebreak\tiny
\be\label{*} \tabcolsep=3.5pt
 \begin{tabular}{|l|l|}
 \hlx{hv}
 $\per\in\Aut \fosp_{\eps_3} (4|2), \text{ where }
\eps_3:= \fnfrac{1}{2}(-1\pm i\sqrt{3})$&
 $\per(a,u)=((a_3,a_1,a_2), \eps  u_3\otimes u_1\otimes
 u_2)$\\
 \hlx{vhv}
  $d_{23}\in\Aut \fosp_\lambda(4|2),\text{ where }\RE\lambda=-\frac12$&
 $d_{23}(a,u)=((a_1,a_3,a_2), u_1\otimes u_3\otimes u_2)$ for any \\[1pt]
&
 $\arraycolsep=0pt\begin{array}{l}
(a_1,a_2,a_3)\in\fsl(V_1)\oplus\fsl(V_2)\oplus\fsl(V_3)= \fosp_\alpha(4|2)_{\ev}\\
 \text{and}\ u_1\otimes u_2\otimes u_3\in V_1\otimes V_2\otimes V_3=\fosp_\alpha(4|2)_{\od}
\end{array}$
 \\
 \hlx{vhv}
 $A\in\Aut\fpo(0|2n)$, \text{~~ where~~}$\fpo(0|2n)\simeq\Cee[\theta_1, \dots, \theta_{2n}]$&
 $A(\theta_i)=(-1)^{\delta_{1i}}\theta_i$,\\

 $B\in\Aut\fpo(0|2n)$&
$B(\theta_i)=\theta_i + \partial_{\theta_i}(\theta_1\ldots\theta_{2n})$,\\  
 \hlx{vhv}
 $\delta_\lambda\in \begin{cases} \Aut\fvect(0|n), \text{~~ where~~
 }\fvect(0|n)=\fder\Cee[\theta_1, \dots, \theta_n]\\ \Aut\fgl(n|m),\quad \text{where
 }\lambda\in\Cee^\times\end{cases}$ &
 $\arraycolsep=0pt
 \begin{array}{l}
 \delta_\lambda(\theta_i)=\lambda\theta_i\text{~~ for all~~
 }i\\
 \delta_\lambda\begin{pmatrix}A&B\\C&D\end{pmatrix}=\begin{pmatrix}A&\lambda B\\\lambda^{-1}C&D\end{pmatrix}\end{array}$\\
 \hlx{vhv}
 $\Pty\in  \Aut\fg$ &
 $\Pty(x)=(-1)^{p(x)}$ for any $x\in\fg$\\
\hlx{vhv}
 $\Ad_{J_{k,2n}(A)}\in \Aut\fosp(k|2n) \text{ where }  J_{k,2n}(A)=\diag(A,1_{2n})$&
 $\arraycolsep=0pt
 \begin{array}{l}
\text{for any } A\in\text{O}(k) 
\text{ such that }   \det A=-1,\  AA^t=1
 \end{array}$\\
 \hlx{vhv}
 $-\st\in\Aut\fgl(n|m)$&
 $-\st\begin{pmatrix}A&B\\C&D\end{pmatrix}=\begin{pmatrix}-A^t&C^t\\ -B^t&-D^t\end{pmatrix}$\\
  \hlx{vhv}
  $\Pi\in\Aut\fgl(n|n)$& $\Pi\begin{pmatrix}A&B\\ C&D\end{pmatrix}=\begin{pmatrix}D&C\\ B&A\end{pmatrix}$\\
 \hlx{vhv}
 $q\in\Aut\fq(n)$& $q: (A, B)\longmapsto(-A^t, i B^t)$\\
 \hlx{vh}
\end{tabular}
\ee
}
\end{table}

{\MathSkip{.5}%
\begin{table}[ht]\centering
{\footnotesize%
\caption{Lie superalgebras $\fg^{\ell(m)}_\psi$ and components
$\fg_i$ for $0\leq \bar i\leq m-1$ (after \cite{Se,LSS})} \label{t-add-3.1}
\begin{minipage}{\textwidth}
\hskip 1.5em The~$\fosp(m|2n)=\fosp_B(V)$-module $S^2(V)$ is reducible: $S^2(V)=S_0^2(V)\oplus \Cee B$. 
Set $T_{m,n}:=\diag(-1, 1_{2m-1}, 1_{2n})$ for $m>0$; let $J_{2n}=\antidiag(1_n, -1_n)$, set $B_{m,2n}:=\diag(1_m, J_{2n})$. The components $\fg_{\bar i}$ are described as $\fg_{\bar 0}$-modules. Let $\eps_3:= \fnfrac{1}{2}(-1\pm i\sqrt{3})$. In line $8$, let $L^k$ be the irreducible $\fosp(1|2)$-module with highest weight $k$ and even highest weight vector. For $\fg_{\bar 0}\subset\fgl(V)$ we denote by $\id$ the tautological $\fg_\ev$-module $V$; let $\ad$ be the adjoint $\fg_{\bar 0}$-module. In line 1, $\fg$ is any simple finite-dimensional Lie superalgebra.
\end{minipage}

\vskip2mm
\be\label{**}
 \tabcolsep=3.1pt
\begin{tabular}{|c|c|c|c|c|c|c|}
 \hlx{hv}
& $\fg$ & $\psi$ & $\fg_{\bar 0}$ & $\fg_{\bar 1}$ & $\fg_{\bar 2}$ & $\fg_{\bar 3}$\\
 \hlx{vhv}
$1$& $\fg$ & $\id$ & $\fg$ & $-$ & $-$ & $-$\\
 \hlx{vhhv}
$2$&  $\fsl(2m|2n)$ & $\Ad_{B_{2m,2n}}\circ(-\st)$ & $\fosp(2m|2n)$ &
 $S_0^2(\id)$ & $-$ & $-$\\
 \hlx{vhv}
$3$& $\fsl(2m+\mkern-2mu1|2n)$ & $\Ad_{B_{2m+\mkern-2mu1,2n}}\circ(-\st)$
 & $\fosp(2m+\mkern-2mu1|2n)$ & $S_0^2(\id)$&$\Pi(\id)$   & $\Pi(\id)$\\
 \hlx{vhv}
$4$&  $\fpsl(2n|2n)$ & $\Ad_{B_{2n,2n}}\circ(-\st)$ & $\fosp(2n|2n)$
 & $S_0^2(\id)$ & $-$ & $-$\\
 \hlx{vhv}
$5$&  $\fpsl(n|n)\text{  for  } n>2$ & $\Pi$ & $\fp\fq(n)$ & $\ad^*$ & $-$ & $-$\\
 \hlx{vhv}
$6$&  $\fpsl(n|n)\text{  for  } n> 2$ & $\Pi\circ(-\st)$ & $\fspe(n)$ &
 $\ad^*$ & $-$ & $-$\\
 \hlx{vhv}
$7$&  $\fosp(2m|2n)$ & $\Ad_{T_{m,n}}$ & $\fosp(2m-\mkern-2mu1|2n)$ & $\id$ & $-$ & $-$\\
 \hlx{vhv}
$8$&  $\fosp_{\eps_3}(4|2)$ & $\per$ & $\fosp(1|2)$ &$\ad=L^2$
 &$\Pi(L^3)$& $-$\\
 \hlx{vhv}
 $9$&  $\fh^{(1)}(2n)\text{  for  } n>2$ & induced by $A\in\Aut\fpo(0|2n)$ & $\fh(2n-\mkern-2mu1)$ & $\ad^*$ & $-$ & $-$\\
\hlx{vhhv}
 $10$&  $\fpsq(n) \text{  for  } n>2$ & $q$ & $\fo(n)$ & $\Pi(S_0^2(\id))$ &
 $S_0^2(\id)$ & $\Pi(E^2(\id))$\\
 \hlx{vhv}
$11$&  $\fpsq(n) \text{  for  } n>2$ & $\delta_{-\mkern-2mu1}$ & $\fsl(n)$ & $\Pi(\ad)$ & $-$ & $-$\\
\hlx{vh}
\end{tabular}
\ee} \vspace{-3mm}
\end{table}
}

\normalsize

\vfill
\newpage



\end{document}